\documentclass[12pt,draft]{article}
\usepackage{amsfonts}
\usepackage{amssymb}
\usepackage{amsmath}
\usepackage{amsthm}

\numberwithin{equation}{section}
\newtheorem{lemma}{Lemma}[subsection]
\newtheorem{theorem}[lemma]{Theorem}
\newtheorem*{theorem*}{Theorem}
\newtheorem{proposition}[lemma]{Proposition}
\newtheorem{corollary}[lemma]{Corollary}
\theoremstyle{definition}
\newtheorem{definition}[lemma]{Definition}
\theoremstyle{remark}
\newtheorem{remark}[lemma]{Remark}
\newtheorem{example}[lemma]{Example}

\newcommand{\Fr}{Fr\'echet }
\newcommand{\Hor}{H\"ormander}
\newcommand{\Ess}{\text{Ess}}

\newcommand{\aPsDO}{{\Psi}} 
\newcommand{\cPsDO}{{\it\Psi}} 

\newcommand{\norm}[1]{\lVert#1\rVert}
\newcommand{\abs}[1]{\left|#1\right|}
\newcommand{\pair}[1]{\left\langle#1\right\rangle}

\newcommand{\R}{\mathbb{R}}
\newcommand{\C}{\mathbb{C}}
\newcommand{\Z}{\mathbb{Z}}
\newcommand{\Q}{\mathbb Q}
\newcommand{\N}{\mathbb N}
\newcommand{\A}{\mathcal{A}}

\newcommand{\G}{\mathcal{G}}
\newcommand{\K}{\mathcal{K}}

\newcommand{\sS}{\mathcal{S}}
\newcommand{\sN}{\mathcal{N}}
\newcommand{\sC}{\mathcal{C}}

\newcommand{\sH}{\mathcal{H}}
\newcommand{\sK}{\mathcal{K}}
\newcommand{\sL}{\mathcal{L}}

\newcommand{\g}{\mathfrak{g}}

\newcommand{\pd}{\partial}
\newcommand{\x}{\times}
\newcommand{\rx}{\rtimes}

\newcommand{\sm}{\backslash} 
\newcommand{\mz}{\backslash\{0\}} 
\newcommand{\ep}{\varnothing}

\newcommand{\twobytwo}[4]
{\left[\begin{matrix} #1 & #2\\ #3 & #4\end{matrix}\right]}

\begin{document}

\title{Transversally Elliptic Operators}
\author{Xiaodong Hu}
\date{May, 2004}
\maketitle

\begin{abstract}

In this paper we investigate the index theory of transversally elliptic
pseudo-differential operators in the framework of noncommutative geometry.

We give examples of spectral triples in the sense
of Alain ~Connes and Henri ~Moscovici in \cite{CM:local}
that are transversally elliptic but non-elliptic.
We prove that these spectral triples
satisfy the conditions which ensure the Connes-Moscovici
local index formula applies.
We show that such a spectral triple
has discrete dimensional spectrum.

We introduce an algebra $\cPsDO^\infty\rx G$ consisting
of families of of pseudo-differential operators,
but with convolution like product over parameter space $G$.
We show the symbolic calculus of this algebra,
which is similar to $\cPsDO^\infty$, the algebra of 
pseudo-differential operators.

Given a spectral triple on $\cPsDO^\infty\rx G$,
we show that there is a finite number of trace-like functionals
$\tau_0, \tau_1, \ldots, \tau_N$ on $\cPsDO^\infty\rx G$
which are defined and used in \cite{CM:local}
for the computation of the local index formula.
Those $\tau$-functionals generalize the Wodzicki
residue of $\cPsDO^\infty=\cPsDO^\infty\rx \{1\}$,
which is $\tau_0$, the only nonzero among these $\tau$ functionals.
Only the last nonzero $\tau_N$ is a trace on $\cPsDO^\infty\rx G$.
It is shown $N$ is bounded by the sum of the dimensions of the
compact Lie group $G$, and the underlying manifold $M$ on
which $G$ acts. Moreover those $\tau$-functionals,
evaluated at $A\in\cPsDO^\infty\rx G$,
are determined by the transversal symbol of $A$.

The calculus on $\cPsDO^\infty\rx G$ is helpful
in the computation of Connes' Chern character of
a spectral triple. We show that Connes-Moscovici
local index formula still works in the transversally
elliptic case.

\end{abstract}

\section{Introduction}

\subsection{Background}

In this section we briefly review some historical
background that motivated and influenced our work the
most.

Atiyah introduced in a 1974 lecture note
\cite{Atiyah:TEO} the index of an invariant transversally
elliptic pseudo-differential operator $P$
relative to a compact Lie group
$G$ action on a compact manifold $M$.
This index generalizes the character
of the index of an invariant elliptic operator,
which is a smooth, central function on $G$.
Although the kernel and co-kernel
of a transversally elliptic operator $P$ are no longer
finite dimensional as in the elliptic case,
the index
\begin{equation*}
\label{index}
index^G(\sigma(P)) =ker(P)-ker(P^*)\in C^{-\infty}(G)^{G}
\end{equation*}
makes sense as a central distribution on $G$.
The index depends on the equivariant $K$-theory class
of the transversal part of its principal symbol.
An explicit index formula was given
in \cite{Atiyah:TEO} for torus acting with finite isotropy.

The operator algebra approach to the
transversally elliptic operators
followed soon afterwards.
P. Julg (1982, \cite{Julg})
gave the following observation.
Just like an elliptic pseudo-differential operator
gives a $K$-homology class in $KK(C(M),\C)$,
a transversally elliptic pseudo-differential operator
naturally induces a $K$-homology class in $KK(A,\C)$, where
the algebra $A$ is the crossed product algebra $C(M)\rx G$.

In his seminal 1985 paper \cite{Co:NDG},
where noncommutative differential geometry
was introduced,  Alain Connes used
smooth cross product algebras as important
examples for noncommutative geometry.
Transversally elliptic operators for
foliations were studied as important 
examples of Fredholm modules over
smooth groupoid algebras. The index
of transversally elliptic operator
problem evolves into the Connes-Chern
character problem for summable Fredholm
modules. The group action case and the foliation
case are not identical although quite similar.
The common ground between the two is
when the group action is free.
So it became natural to study 
the smooth crossed product algebra 
$\A=C^\infty(M) \rx G$.
Both the foliation smooth algebra
and the group action smooth algebra
can be viewed as smooth groupoid
algebras. 

Smagin and Shubin (1987, \cite{Smagin-Shubin})
introduced the formal zeta function constructed from
the transversally elliptic symbols. They showed the use
of the resolution of singularities of the phase function.
Later Yuri Kordyukov has studied
Fredholm modules for group action algebra
and spectral triples for Riemannian foliation algebra
(see, \cite{Kor:GTEO1}, and \cite{Kor:GTEO2}).

Connes and Moscovici (1995, \cite{CM:local})
gave the general local index formula
for spectral triples, together with
a new working example for hypo-elliptic operator.
Later they showed many important applications.
We intend to use this result to give an
index formula for transversally elliptic
operators on a smooth manifold
relative to a compact Lie group action.

Also, we have been inspired also by some
work done essentially in classical
(instead of operator algebraic) approach.
For instance, Helga Baum (1983, \cite{Baum}) studied
the transversal index of pseudo-Riemannian
Dirac operators; also Berline and Vergne (1997, \cite{BV})
computed Chern characters of so called
transversally elliptic good symbols.

\subsection{The main results}

\subsubsection{The spectral triple}

The noncommutative ``space'' that we consider, is the
smooth crossed product algebra $\A=C^\infty(M\rx G)$,
which is a \Fr algebra with a natural topology.
The elements of $\A$ are smooth functions
(denoted by $a$, $b$) on $M\x G$,
with the product
\begin{equation*}
  (a*b)(x,g)=\int_G a(x,h)b(h^{-1}x, h^{-1}g)d\mu(h).
\end{equation*}
$\A$ can be viewed as the smooth convolution
algebra of the groupoid $M\rx G$ induced
by the group action, which 
is the replacement for (the algebra of functions on)
the quotient space $M/G$.
This algebra has been shown to be closed under
holomorphic functional calculus (see \cite{Schweitzer}).

Let $D$ be a first order
transversally elliptic pseudo-differential operator.
$\sH$ is usually the graded direct sum of $L^2$-sections
of the complex vector bundles on which $D$ acts
as an unbounded self-adjoint operator.
Let $\epsilon$ be the grading operator of $\sH$.
Further we assume that $D^2$ has scalar symbol.
Although this assumption is rather strong, there
are always examples of such operators, such as the
Dirac operator associated to a general Clifford
module and a Clifford connection. 
$\A$ acts upon $\sH$ by extending the
$G$-action $\rho$ on the fibers
of the complex bundles:
\begin{equation*}
 (\rho(a)\cdot s)(x)=\int_G a(x,g)(\rho(g)s)(g^{-1}x)dg.
\end{equation*}
Thus for any $a\in\A$, $\rho(a)\in B(\sH)$,
$\rho(a)\epsilon=\epsilon\rho(a)$, and
$D\epsilon=-\epsilon D$.

A pseudo-differential operator $K$ on $\sH$
is called transversally smoothing
(to be denoted by $K\in \K_\A$)
if for any $a\in \A$, $\rho(a)K$ is a smoothing
operator (in particular trace class).

We first show that $(\A,\sH,D)$ is a spectral triple
in the following sense:
\newcounter{tmp}\setcounter{tmp}{1}
(\roman{tmp}) for all $a\in\A$, $[D, \rho(a)]\in B(\sH)$;
\stepcounter{tmp}
(\roman{tmp}) for any $a\in\A$, $\rho(a)(1+\abs{D})^{-1}$
is compact;
\stepcounter{tmp}
(\roman{tmp})  there is a transversally smoothing operator
$K\in\K_\A$ such that $\abs{D}+K$ is invertible, and
the inverse  $(\abs{D}+K)^{-1}$ is in the class
$\sL^{(p, \infty)}$.

Recall that for $p=\dim{M}>1$, 
where $\sL^{(p, \infty)}$ is the ideal
of $B(\sH)$ consisting of those
compact operators $T$ whose $n$-$th$
characteristic value
satisfies
\begin{equation}\label{eq:ppsummable}
  \mu_n(\abs{T}) = O(n^{-1/p}).  
\end{equation}
When $p=1$, although (\ref{eq:ppsummable}) is stronger than
$T\in \sL^{(1,\infty)}$, it is also satisfied.

\subsubsection{Noncommutative residues}
\label{sec112}
Throughout this and the next section,
we assume we have fixed a spectral triple
$(\A, \sH, D)$ as described in the previous
section.

From the algebra $\A$ and its representation on $\sH$,
we introduce the crossed product algebra
of pseudo-differential operator and the group $G$:
\[
\cPsDO_G=\cPsDO^\infty(E)\rx G =
\bigcup_k \cPsDO^k(E)\rx G.
\]
The operators in this algebra are families of
pseudo-differential operators with parameter
space $G$ and convolution-like composition maps.
We show that, similar to ordinary
pseudo-differential operators, on $\cPsDO_G$
there is a similar filtration
(in terms of order $k$) and symbolic calculus.

There is a distribution kernel for any operator in $\cPsDO_G$.
The wave front set of such a kernel has a certain form
(see \ref{lemma:xpsdo}), 
just like the wave front set of a pseudo-differential operator,
which has in a diagonal form (the micro-local property).

The motivation for introducing these operators is
the computation of the Connes-Moscovici local
index formula, and is described in section \ref{sec113}.

For the transversally elliptic $D$ in the spectral
triple, $\abs{D}$ is only transversally elliptic.
But as mentioned, we may add an $G$-invariant,
transversally smoothing $K$ to $\abs{D}$
so that the replacement $\abs{D}+K$ is invertible and
elliptic. By simple wave front set computation,
it is clear that $K$ does not affect
the asymptotic trace integral formulas.
So for the convenience of argument we may ignore $K$.

On the algebra $\cPsDO_G$, we study the
noncommutative residues defined in
Connes-Moscovici \cite{CM:local}.
For $P\in\cPsDO^k(E)\rx G$, the zeta function
\begin{equation}\label{intro:zeta:trace}
\zeta_{P,D}(z)=Trace( P \abs{D}^{-2z})
\end{equation}
is initially defined and analytic on
the half plane
$\{z\in \C : 2 Re(z)>k+\dim{M}\}$.

For our spectral triple,
we define the {\it dimension spectrum relative to $\cPsDO_G$}
as the minimal closed subset of $\C$ on the
complement of which $\zeta_{P,D}$ can be extended
to a holomorphic function for any $P\in \cPsDO_G$.
$\cPsDO_G$ contains $\A_D$ as defined by
Connes and Moscovici.
The dimension spectrum could be larger if we replace
$\cPsDO_G$ by one of its subalgebra.

Analysis of behavior the of the
trace (\ref{intro:zeta:trace}) amounts
to some oscillatory integrals with a
phase function related to the fixed point
sets of the $G$ action on $M$.

\begin{theorem}
The dimension spectrum relative to $\cPsDO_G$
is a discrete subset of the nationals $\Q$.
\end{theorem}

For any $P\in\cPsDO^\infty(E)\rx G$ and $q=0,1,\ldots$
the noncommutative residues $\tau_q$, $q=0,1,\ldots$
are defined as

\begin{equation}
\tau_q(P)=\tau^{\abs{D}}_q(P)=Res_{z=0} (z^q \zeta_{P,D}),
\end{equation}
that is, the residues of $z^q \zeta_{P,D}$ at $z=0$.
These residues may depend on the choice of $D$.
For convenience, $\tau_{-1}(P)$ is also defined this
way, it is the value of the zeta function at zero
if zero is not a pole.

\begin{theorem}
For any $P\in \cPsDO_G$, the poles of $\zeta_{P,D}$
are of multiplicity at most $\dim{M}+\dim{G}$.
Therefore, there are only finitely many nonzero
noncommutative residues $\tau_q$, $q=0,1,\ldots,$
up to possibly $\dim{M}+\dim{G}-1$.
\end{theorem}

The dimension spectrum and the maximal $q$ (for
$\tau_q$ to be nonzero) depends primarily on
the action of $G$. 
It is well known that when $G$ is the trivial group, where
$\cPsDO_G=\cPsDO$ is the algebra of pseudo-differential operators,
only $\tau_0$ is nonzero and it is the Wodzicki residue.
The Wodzicki residue extends any one of the Dixmier traces on
the subalgebra of pseudo-differential operators of order $-\dim{M}$ to
the algebra of all pseudo-differential operators.

It is shown in \cite{CM:local} that only the last $\tau_q$ is a trace.
As shown in \cite{CM:local},
when $q>0$, $\tau_q$ always vanishes on the Dixmier trace
class operators (that is, $\sL^{(1,\infty)}$) and
in particular, on $\cPsDO^{-\dim{M}}\subset \sL^{(1,\infty)}$.
But on $\cPsDO_G^{-\dim{M}}$, $\tau_0(P)$ is equal to any Dixmier trace.
Therefore $\tau_0$ still extends any Dixmier trace  as a functional but
it may not be a trace on $\cPsDO_G$.

\begin{theorem} 
For any $q=0,1,\ldots,\dim{M}+\dim{G}-1$,
$\tau_q$ vanishes on $\sK_G$, the ideal of
transversally smoothing operators.
\end{theorem}

\begin{theorem}
For any $q=0,1,\ldots,\dim{M}+\dim{G}-1$,
and any $P\in \cPsDO_G$, $\tau_q$ depends only
on the transversal part of the full symbols of $P$
of order no lower than $-\dim{M}$.
\end{theorem}

The above theorem asserts that the noncommutative residues
are in principle ``computable'' in terms of symbolic calculus.

\subsubsection{Connes-Moscovici local index formula}
\label{sec113}

The index formula of $(\A,\sH,D)$
follows from the Connes-Chern character
formula in periodic cyclic cohomology.
The operator-theoretical version of the local index
formula is given in Connes-Moscovici \cite{CM:local}
in full generality.
With the the theorems in section \ref{sec112}
we show that it is possible to apply
the Connes-Moscovici local index
formula for our specific spectral triple
$(\A, \sH, D)$.

For any operators $A\in \A$, we use the
following notations
\begin{equation*}
dA=[D,A],\;
\nabla(A)=[D^2,A],\;
A^{(k)}=\nabla^k(A).
\end{equation*}
Then all the above operators are in $\cPsDO_G=\cPsDO^\infty(E)\rx G$.
In particular, the operators
\begin{equation}\label{eq1:1:3}
a^0 (da^1)^{(k_1)} \ldots (da^n)^{(k_n)}
\end{equation}
are in $\cPsDO^\infty(E)\rx G$,
where $a^0,\ldots, a^n \in \A$,
acting on $\sH$ by $\rho$.
Let $\A_D$ be the subspace
of $\cPsDO^\infty(E)\rx G$ generated
by those operators in (\ref{eq1:1:3}).

$\abs{D}^{-1}$ (which is really $(D+K)^{-1}$),
is in $\cPsDO^{-1}(E)\rx G$ and we have
\begin{equation}\label{eq1:1:4}
  a^0 (da^1)^{(k_1)} \ldots (da^n)^{(k_n)} 
  \abs{D}^{-2\abs{k}-n} \in \cPsDO^0(E)\rx G.
\end{equation}
In particular, the left hand side of (\ref{eq1:1:4}) is bounded.
So we assert the Connes-Moscovici local index formula
holds.

\begin{theorem}
Let $(\A, \sH, D)$ be an even spectral triple 
defined by a first order transversally elliptic
pseudo-differential operator $D$ and with
all the above conditions.

The Connes character $ch(\A,\sH,D)$
in periodic cyclic cohomology is represented by the following
even cocycle in the periodic cyclic cohomology:
\begin{equation}
\begin{split}
\phi_{2m}(a_0,\ldots, a_{2m})
&=\sum\limits_{k\in\Z^{2m},q\ge 0}
c_{2m,k,q}\cdot\\
&\tau_q \left( a^0(da_1)^{(k_1)}\cdots (da_{2m})^{(k_{2m})}
\abs{D}^{-2\abs{k}-2m}\right)
\end{split}
\end{equation}
for $m>0$ and
\begin{equation}
\phi_0(a^0)=\tau_{-1}(\gamma a^0).
\end{equation}
In the above formula
$k=(k_1,\ldots, k_{2m})$ are multi-indices
and $c_{2m,k,q}$ are universal constants
given by
\begin{equation}
c_{2m,k,q}=
\frac{(-1)^{\abs{k}}}{k!\tilde{k}!}\sigma_q(\abs{k}+m),
\end{equation}
where $k!=k_1!\ldots k_{2m}!$,
$\tilde{k}!=(k_1+1)(k_1+k_2+2)\ldots (k_1+\ldots+k_{2m}+2m)$,
and for any $N\in\N$, 
$\sigma_q(N)$ is the $q$-$th$ elementary polynomial
of the set $\{1,2,\ldots,N-1\}$.
\end{theorem}

There is a similar theorem for odd spectral triples. 

\subsection{Outline of the paper}

In section 2 we review mostly known facts
we need, including wave front set
and pseudo-differential operators.
Simple wave front set calculation
allows us to use the spectral analysis of
elliptic operators instead of
transversally elliptic ones.
In section 3 we first review Atiyah's definition
of a transversally elliptic operator
operators. Then we define the crossed product
algebra $\A$, our noncommutative space.
Next the algebra of crossed products with
pseudo-differential operators is introduced.
In the rest of this section we establish
basic rules of calculus as needed later.
The notion of transversally smoothing operators
is introduced to simplify arguments.
In section 4 we discuss trace formulas.
First we study the trace class operators
and the trace formula in terms of the
distributional kernel. Then we show
how the asymptotic expansion of oscillatory
integrals comes into the picture. Next we
review some known theorems about the asymptotic
expansion of oscillatory integrals.
The definition noncommutative residue is reviewed
and then their properties are discussed.
This allows us to develop our main statements
about the noncommutative residues.
Further techniques are needed to decide
these residues by the symbol.
In section 5 we discuss
the $K$-theoretic aspect of the index of
the transversally elliptic operators. Here we show that
the index of a transversally elliptic operator
as defined by Atiyah \cite{Atiyah:TEO}
is a natural part of the $K$-theory of the
noncommutative space $\A$.
Next we show that the Connes-Chern character in $K$-theory
is related to the transversal index.
Finally we apply our residue formula to
the Connes-Moscovici local index formula for
the spectral triples decided by transversally
elliptic operators.

\subsection{Acknowledgments}

The author is deeply in debt to Professor Henri Moscovici
for his kind instructions, discussions, encouragement.
The author also wishes to thank Professor Dan Burghelea,
Professor Alexander Dynin, Professor George Elliott, and
Dr. Hanfeng Li for very helpful discussions, critics
and enhancements.

\section{Preliminaries}

\subsection{Conventions}

Unless otherwise specified, we use the following notations
throughout the paper.

Let $G$ be a compact Lie group acting
smoothly on a compact smooth manifold $M$.
We denote by $\g$ the Lie algebra of $G$.
For convenience, we assume that $M$ is a compact, oriented
Riemannian manifold with an invariant
Riemannian metric and $G$ acts by orientation
preserving isometries.
We also assume that a fixed $G$-invariant volume
form on $M$ is selected. For an element $g\in G$,
let $M^g=M(g)$ be the fixed point set of $g$.
(which under all the above assumptions is a disjoint union of
totally geodesic submanifolds of $M$).
We will consider complex hermitian $G$-equivariant
bundles over $M$, which are assumed here to have 
smooth hermitian metrics preserved by the $G$-action.

Let $\N$ be the set of the natural numbers, including zero.

For a smooth manifold, let $C^\infty_c(M)$
be the compact supported, complex valued test functions
with the induced limit topology of local test functions
(often denoted by ${\mathcal D}(M)$),
and $C^{-\infty}(M)$  the corresponding
distributions(often denoted by ${\mathcal D}'(M)$).
We work mostly on compact manifolds,
though. For a hermitian bundle $E$ over $M$,
$\Gamma^\infty_c(E)$ (or ${\mathcal D}(M, E)$) and
$\Gamma^{-\infty}(E)$ (or ${\mathcal D}'(M, E)$)
the sections of test functions and distributional
sections.

For general notations or definitions used
in this paper, the standard reference is
Connes' book \cite{Co:book}. We refer to
papers \cite{AB:Lef1} and \cite{AS1} and
books \cite{Ho1} for other definitions
and backgrounds.

\subsection{Wave front set}

\subsubsection{Definition}
For a distribution $u$ with compact support
to equal a smooth function, its Fourier transform
needs to satisfy: for any $N\in\N$ there
is a $C_N$ such that
\begin{equation}\label{eq:local:not:wf}
\abs{\hat{u}(\xi)} \le C_N(1+\abs{\xi})^{-N} 
\end{equation}
for all $\xi\in \R^n$.

\begin{definition}
Let $u$ be a distribution on $\R^n$ with compact
support. Then at $x\in \R^n$ its wave front set $\Sigma_x$
is a cone of all $\eta\in \R^n\mz$ having
no conic neighborhood $V$ such that (\ref{eq:local:not:wf})
is true for all $\xi \in V$.
For an arbitrary distribution $u$ on an open subset $X$ of
$\R^n$ let
\begin{equation}
\Sigma_x(u)=\bigcap_{\phi\in C^\infty_c(\R^n),\phi(x)\not=0}
\Sigma_x(\phi u).
\end{equation}

The wave front set for $u$ is
\begin{equation}
WF(u)=\{(x,\xi)\in X\x (\R^n\mz); \xi\in \Sigma_x(u)\}.
\end{equation}
\end{definition}

For a smooth vector bundle $E$ over $M$,
let $\Omega(M)$ be the volume bundle over $M$
and $E'=Hom(E, \Omega(M))$.
By definition, a distributional section
is an element of $\Gamma^{-\infty}(M,E)$, the
dual space of $\Gamma^\infty_c(M,E')$.

The invariant definition for vector-valued
distributions is given below. For details,
see, for example, \cite{Ho1},
\cite{Duistermaat:FIO} or \cite{BEP}.

\begin{definition}
For a smooth bundle $E$ on $M$,
let $u\in \Gamma^{-\infty}(M,E)$ be a
distributional section of $E$.
The wave front set $WF(u)$ is
a subset of $T^*M\mz$, such that
for any $(x,\xi)\in T^*M\mz$, $\xi_x \not\in WF(u)$
if and only if for any (phase function)
$\psi\in C^\infty_c(M\x \R^p, \R)$,
$d\psi(.,y)_x=\xi_x$, there is
an $s\in \Gamma^\infty_c(M, E')$, $s(x)\not=0$,
and a neighborhood $U_y$ of $y$ in $\R^p$,
such that for all $n\in \N$, and all $y\in U_y$,
\begin{equation}
\abs{\langle u, exp(-it\psi(x,y'))s\rangle}
=O(t^{-n})
\end{equation}
when $t \to \infty$.
\end{definition}

Wave front set is a refinement of the
singular support for distributions.

\begin{proposition}
$WF(u)$ is a closed conic subset.
The projection of $WF(u)$ onto $M$
is the singular support of $u$.
\end{proposition}

\subsubsection{Push-forward of a distribution and its wave front set}
For distributions, the push-forward
is a natural operation, and the wave front
set under this operation is described below.

\begin{definition}\label{def:push}
Let $f: M \to N$ be a proper smooth map,
Let $S$ be a conic subset of $T^*M\mz$,
then the push-forward $f_*S$ of $S$
is
\begin{equation}
f_*S=\{ (y, \eta)\in T^*N\mz:
  \exists x \in M, y=f(x),
  (x, f^*(\eta))\in (S\cup 0_{T^*M})
\}.
\end{equation}
\end{definition}

\begin{theorem}
Let $f: M \to N$ be a smooth map,
$u\in C^{-\infty}(M)$ such that
$f|_{supp(u)}$ is proper. Then
the push-forward $u\mapsto f_*(u)$
is a well-defined continuous linear map
\begin{equation}
f_*: C^{-\infty}(M) \to C^{-\infty}(N).
\end{equation} 
Moreover,
\begin{equation}
WF(f_*u)\subset f_*(WF(u)).
\end{equation}
\end{theorem}

\begin{example}
(1) Take $M=\{p\}$, the $0$-dimensional
manifold of a single point, $f$ is simply
decided by its image $y=f(p)$ in $N$.
Let $1_M$ be the constant
function on $M$ with value $1$. Then $f_*(1_M)$
is simply $\delta_y$, the delta function
of $N$ at the point $y$.
Simple computations right from definitions show
that
$WF(\delta_y)=T^*_y N\mz$,
$WF(1_M)=\ep$,
and $f_*(\ep)=T^*_y N\mz$.

(2) A little more generally, let $f$ be an embedding,
let $\alpha$ be a compactly supported
smooth function on $M$,
$f_*(\alpha)=\alpha\mu_M$ is usually
called the density of $M$ in $N$,
direct computation shows
\begin{equation}
  WF(f_*\alpha)=\{(x,\xi) \in T^*N,
  x\in supp(\alpha), \xi \in \sN^*_M(N)\mz\},
\end{equation}
where  $\sN^*_M(N)=ker(f^*)$ is the conormal bundle
of $M$ in $T^*N$. Compare $WF(\alpha)=\ep$
and $f_*(\ep)=ker(f^*)$.

(3) Let $\pi: M\x N \to N$ be the projection to $N$.
Let $M$ be compact so the push-forward of $\pi_*$
is well-defined.
For $v\in C^\infty_c(N)$ and $u\in C^{-\infty}(M\x N)$,
\begin{equation}
(\pi_*u,v)=(u, \pi^*v)
=(u, 1_M\otimes v)
=\int_{M}(u(x,\cdot), v(\cdot))dx, 
\end{equation}
the last formal integral is true when $u$ is smooth.
$WF(\pi_*u)$ is the projection of $WF(u)$ onto $T^*N$.

In fact, the push-forward operation contains the
integration for a distribution as a special case.
Let $N=\{pt\}$ be a single point
manifold so the functions and distributions on it are
just complex numbers. Let $c$ be the unique map $M\to N$,
$c_*(u)$ is defined when $M$ is compact and it is
just the number
\begin{equation}
c_*(u)=(u, 1_M)=\int_M u(x)dx.
\end{equation}
where the last integral is formal.

(4)
Now let $f: M \to N$ be a submersion from a compact
manifold $M$ to compact manifold $N$.
Let $u$ be a smooth function on $M$, so $WF(u)=\ep$.
Then the push-forward $f_*u$ is still an integration along the
vertical fibers with our chosen volume form.
\begin{equation}
WF(f_*u)=f_*(WF(u))=\ep,
\end{equation}
so as we know $f_*u$ is smooth. When $u$ is a distribution
$f_*(WF(u))$ is again the projection along vertical
directions.
\end{example}

It is straightforward to extend the above theorem to
maps between bundles, because the wave front set
is a local concept.

\begin{theorem}\label{thm:push}
Let $f: (E, M) \to (F, N)$
be a smooth bundle map between two
hermitian bundles,
$u\in \Gamma^{-\infty}(E)$ such that
$f|_{supp(u)}$ is proper. Then
the push-forward $u\mapsto f_*(u)$
is a well-defined continuous linear map
\begin{equation}
f_*: \Gamma^{-\infty}(E,M) \to \Gamma^{-\infty}(F,N).
\end{equation} 
Moreover,
\begin{equation}
WF(f_*u)\subset f_*(WF(u)).
\end{equation}
\end{theorem}

\subsubsection{Schwartz kernel and wave front relation}
For two locally convex topological vector spaces $V$ and $W$,
we denote by
\begin{equation}
  \sL_S(U,V)
\end{equation}
the topological vector space of continuous linear maps
from $V$ to $W$, with the strong topology.

Let $(E_1, M_1)$, $(E_2, M_2)$ be
two hermitian bundles, and let
\begin{equation}
A: \Gamma^{-\infty}(E_1) \to \Gamma^{-\infty}(E_2)
\end{equation}
be a linear and continuous operator. 
There is a distributional section
\begin{equation}
  K_A\in \Gamma^{-\infty}(E_2\boxtimes E'_1)
\end{equation}
such that for any $u\in \Gamma^{-\infty}(E_1)$
and $v\in \Gamma^{-\infty}(E'_2)$
\begin{equation}
  \label{eq:x:ker}
  (K_A, v\otimes u) = (v, Au).
\end{equation}
We recall the Schwartz kernel theorem:
\begin{theorem}
  There is a unique continuous linear map
  \begin{equation}
    \sL_S(\Gamma^{-\infty}_c(E_1), \Gamma^{-\infty}(E_2))
    \to
    \Gamma^{-\infty}(E_2\boxtimes E'_1)
  \end{equation}
  that maps $A$ to $K_A$ such that (\ref{eq:x:ker})
  holds. This map is an isomorphism of the 
  topological vector spaces.
\end{theorem}

\begin{definition}
The {\it wave front relation} of $A$ is defined as
\begin{equation}
WF'(A)=\{(\eta_y, \xi_x)\in T^*(M_2\x M_1) |
   (\eta_y, -\xi_x)\in WF(K_A)\}.
\end{equation}
And we define $WF'_{M_i}(A)$,  called $WF'(A)$'s
projection on $T^*M_i$, $i=1,2$, as follows:
\begin{equation}
WF'_{M_1}(A)=\{(x,\xi_x) \in T^*M_1\mz:
  \exists y\in M_2, (y,x,0_y,\xi_x)\in WF'(A)\},
\end{equation}
and
\begin{equation}
WF'_{M_2}(A)=\{(y,\eta_y) \in T^*M_2\mz:
  \exists x\in M_1, (y,x,\eta_y,0_x)\in WF'(A)\}.
\end{equation}
\end{definition}

Note: In the event when $M_1=M_2=M$ (which we encounter
in this paper), to avoid any confusion,
we use an extra label, that is, $WF'_{M, 1}$
and $WF'_{M, 2}$ respectively.

Let $(E_i, M_i)$, $i=1,2,3$ be hermitian bundles
over compact manifolds:
\begin{equation}
A: \Gamma^\infty(E_1) \to \Gamma^{-\infty}(E_2)
\end{equation}
\begin{equation}
B: \Gamma^\infty(E_2) \to \Gamma^{-\infty}(E_3)
\end{equation}
are linear and continuous operators.

\begin{theorem}\label{thm:compo}
If $WF'_{M_2}(A) \cap WF'_{M_2}(B) = \ep$, then
\begin{equation}
B\circ A: \Gamma^\infty(E_1) \to \Gamma^{-\infty}(E_3)
\end{equation}
is a well defined, linear and continuous operator.
Moreover,
\begin{equation}\begin{split}
WF'(B\circ A)
&\subset WF'(B)\circ WF'(A)\\
&\cup (WF'(B)_{M_3}\x 0_{T^*M_1})
\cup (0_{T^*M_3} \x WF'(A)_{M_1}).
\end{split}\end{equation}
(By definition $WF'$ is a binary relation
between $T^*M_2$ and $T^*M_1$ as sets, the
composition ``$\circ$" between wave front
relations is that of the binary relation
between sets.)
\end{theorem}

We will need the following obvious corollary.
Assume now $(E_i, M_i)=(E,M)$ are all the same.

\begin{corollary}
If $WF'_{M, 2}(A) \cap WF'_{M, 1}(B)=\ep$
and $WF'(B)\circ WF'(A) =\ep$, then
$B\circ A$ is a smoothing operator.
\end{corollary}

\begin{theorem}\label{thm:push2}
The wave front relation for the kernel
of $f_*$ in theorem \ref{thm:push} satisfies
\begin{equation}
WF'(f_*)\subset
\{(\xi_x, \eta_{f(x)})\in T^*(M\x N)\mz: 
  f^*(\eta_{f(x)})=-\xi_x \mathrm{\;or \;} f^*(\eta_{f(x)})=0
\}.
\end{equation}
\begin{proof}
We apply theorem \ref{thm:push} to the embedding
$M\to M\x N$, $x\mapsto (x,f(x))$.
\end{proof}
\end{theorem}

\subsection{Pseudo-differential operators}

In this section we review some selected, basic but important
facts about pseudo-differential operators, see \cite{Ho3} for details.
In particularly we need to specify a \Fr topology
on the space of pseudo-differential operators.
This is necessary since we need to integrate over the group $G$ of
a family of pseudo-differential operators, together
with a continuous group action on such families.
In \cite{AB:Lef1} and  \cite{AS1},
a good topology of pseudo-differential operators
are given. Although for most of what we do, the topology on $Op^m$
as in \cite{AS1} is sufficient, the \Fr topology decided by the symbol
is more convenient for our purpose.
Such a \Fr topology has been well known in the early
of symbols and pseudo-differential operators
(see, for example, \cite{MTaylor} and \cite{Cordes}).

\subsubsection{Local calculus}

Recall the definition of scalar symbols
of order $m$ on $\R^n$ (18.1.1 in \cite{Ho3}).
We omit the space in the notations only when it is $\R^n$.
\begin{definition}
For any real number $m$, $S^m=S^m(\R^n\x \R^n)$ is defined
as the set of all $a\in C^\infty(\R^n\x \R^n)$ satisfying
the following condition:
for all multi-indices $\alpha, \beta \in \N^n $ there is
a positive constant $C_{\alpha,\beta}$ such that
\begin{equation}
  p_{\alpha,\beta}(a)=\sup_{x,\xi} \abs{
    \frac{\partial^{\alpha}_\xi\partial^{\beta}_x a(x, \xi)}
    {(1+\abs{\xi})^{m-\abs{\alpha}}}} \le C_{\alpha,\beta}
\end{equation}
for all $x$, $\xi\in \R^n$.
\end{definition}

$S^m$ is a \Fr space with semi-norms $p_{\alpha,\beta}$
for all possible multi-indices $\alpha$ and $\beta$ in $\N^n$.
Each symbol $a\in S^m$
gives rise to a pseudo-differential operator denoted by $Op(a): \sS\to\sS$,
acting on Schwartz functions $\sS$ on $\R^n$: for $u\in\sS$ and
$x\in \R^n$,
\begin{equation}
  Op(a)(u)(x)=(2\pi)^{-n} \int e^{i\pair{x,\xi}}a(x, \xi)
  \hat{u}(\xi) d\xi,
\end{equation}
\noindent $a$ is called the full symbol of $Op(a)$.
This action  extends to a continuous linear operator from
$\sS'$ onto itself. The action of $Op(a)$ is faithful,
that is, $Op(a)=0$ only when $a=0$.

\begin{definition}
We define $\aPsDO^m=\aPsDO^m(\R^n)$ as the topological space of
all operators on $\sS'$ of the form $Op(a)$ for $a\in S^m(\R^n\x \R^n)$,
with topology that of $S^m$. For any $m$,  $\aPsDO^m$ is a \Fr space.
\end{definition}

\begin{theorem}\label{thm:psdo:prod}
If $a_j\in S^{m_j}$, $j=1,2$, then there
exists $b\in S^{m_1+m_2}$ such that
\begin{equation}
Op(a_1)Op(a_2)=Op(b).
\end{equation}
The symbol $b$ is given by
\begin{equation}
  b(x,\xi)
  =e^{i \langle D_y, D_\eta\rangle}
  a_1(x,\eta) a_2(y, \xi) 
  \vert_{\eta=\xi, y=x}
\end{equation}
and it has the following asymptotic expansion:
\begin{equation}
  b(x,\xi) \sim
  \sum_{\alpha}
  \frac{i^{\abs{\alpha}}}{\alpha!}
  D_\xi^\alpha a_1(x,\xi) D_x^\alpha  a_2(x, \xi)
  \vert_{\eta=\xi, y=x} 
\end{equation}
where we recall the notation $D^\alpha_y$
means $(-i)^{\abs{\alpha}}\partial^\alpha_y$.

Moreover, the product map 
\begin{equation}\label{eq:jointclocal}
\begin{split}
\aPsDO^{m_1}\x \aPsDO^{m_2} &\to \aPsDO^{m_1+m_2}\\
Op(a_1)\x Op(a_2) &\mapsto Op(a_1) Op(a_2) 
\end{split}\end{equation} 
is jointly continuous.
\end{theorem}

For the proof, we refer to the proof of 18.1.8 in \cite{Ho3}.
As a consequence, for $m\le 0$, $\aPsDO^m$ is a \Fr algebra.
Obviously, for $m_1<m_2$, $\aPsDO^{m_1}$ embeds naturally into $\aPsDO^{m_2}$.
The union of all $\aPsDO^m$ over all $m\in\R$ is
denoted by $\aPsDO^*$ or $\aPsDO^\infty$, with the induced limit topology,
it is also a \Fr algebra.

\begin{theorem}\label{locinv}
Let $\kappa : X \to X_\kappa$ be a diffeomorphism
between two open subsets of $\R^n$. Then for
any $a\in S^m$ there exists $a_\kappa\in S^m$
such that
\begin{equation}
Op(a_\kappa)=\kappa_* \circ Op(a)\circ \kappa^*
\end{equation}
\end{theorem}

This is the invariance of $\aPsDO^m(\R^n)$ under change
of coordinates (see 18.1.17 of \cite{Ho3}).

\subsubsection{The \Fr Topology}

The notion of pseudo-differential operators extend
to manifolds.

\begin{definition} 
Let $M$ be a smooth manifold, a continuous
operator:
\begin{equation}
P: C^\infty_c(M) \to C^{\infty}(M)
\end{equation}
is called a pseudo-differential operator
($P\in \aPsDO^m(M)$) if for any differentiable
chart of $M$:
\begin{equation*}
\chi: U\subset M \to \chi(U)\subset \R^n
\end{equation*}
and for any $\phi, \psi\in C^\infty(M)$
with supports within $U$, we have
\begin{equation}
P_{\chi, \phi,\psi}=(\chi_*)(\phi P\psi)(\chi^*)\in S^m(\chi(U)).
\end{equation}
\end{definition}

\begin{definition} 
We give $\aPsDO^m(M)$ the topology induced by all
the semi-norms depending on $\phi$ and $\psi$,
\begin{equation*}
\{p_{\alpha, \beta, \phi,\psi}(P) = p_{\alpha,\beta}(\phi P\psi) \}.
\end{equation*}
\end{definition}

\begin{proposition}
For a chosen compatible atlas of charts on a manifold $M$,
let $\phi_i$ be a locally finite partition of unity
subordinate to the covering induced by this atlas.
Then 
\begin{equation*}
\{p_{\alpha, \beta, \phi_i,\phi_j}(P) = p_{\alpha,\beta}(\phi_i P\phi_j) \}.
\end{equation*}
generate the topology of $\aPsDO^m(M)$. Consequently, $\aPsDO^m(M)$
is a \Fr space.
\end{proposition}
\begin{proof}
We need to show that any semi-norm $p_{\alpha, \beta, \phi,\psi}$
is bounded by a finite sum of other semi-norms.
By the partition of unity, 
\[
P=\sum_{i,j} \phi_i P \psi_j.
\]
For any pair $\phi$ and $\psi$
only finitely many of $\phi_i$'s will be needed in the
partition of unity.
This allows us to prove this proposition
with the assumption that $P=\phi_i P \phi_j$.
In other words, we can assume there is one chart
and $P$ has small enough support.

Note that multiplication by $\psi$ is
a pseudodifferential operator of order zero,
that is $\psi\in S^0$, with symbol $\psi(x, \xi)=\psi(x)$.
Using theorem \ref{thm:psdo:prod} for $\phi P$ and $\psi$,
especially the joint continuity, we conclude the
semi-norm $p_{\alpha, \beta, \phi,\psi}$ is bounded by a finite
sum of semi-norms of $\phi P$, times some positive constants
depending only on a finite number of semi-norm of $\psi$ in $S^0$.

Then we repeat this argument with $\phi P$ with $\phi\in S^0$
and $P$. It is routine to check that $\aPsDO^m(M)$ is complete.
With countable semi-norms generating the topology, it is
a \Fr space. 
\end{proof}

\begin{definition}
\begin{equation}
P: \Gamma_c^\infty(M, E) \to \Gamma^\infty(M, F)
\end{equation}
is called pseudo-differential
(to be denoted by $P\in \aPsDO^m(M;E,F)$)
if for a given pair of trivializations $\varphi_E$
of $E$ and $\varphi_F$ of $F$ on $U\subset M$,
and for all $u\in \Gamma_c^\infty(M,E)$, 
there exists $A_{ij}\in \aPsDO^m(U)$ such that
\begin{equation}
\varphi_F(Pu)_i=\sum A_{ij}(\varphi_E u)_j.
\end{equation}
and the topology on $\aPsDO^m(M; E,F)$ over $U$ is given
by the union of semi-norms for $A_{ij}$ over indices $i,j$ on $U$.
\end{definition}

It is straightforward to check that the definition does not
depend on the choice of trivialization of $E$ or $F$;
The topology on $\aPsDO^m(M;E,F)$ is well defined;
and $\aPsDO^m(M; E,F)$ is a \Fr space.

Let $\Omega(M)$ be the volume bundle over $M$.
Recall that, as in \cite{AB:Lef1} and \cite{AS1}, $E'$
is defined as $E^*\otimes \Omega(M)$.
In \cite{Ho3}, \Hor used the bundle $\Omega^{1/2}(M)$ of
half densities to define the adjoint of
a pseudo-differential operator.
With the assumption that $M$ is oriented,
$\Omega$ and hence $\Omega^{1/2}$ are trivial bundles,
so they may be omitted up to an bundle
isomorphism decided by a chosen volume form.
\Hor's approach is equivalent to the one we
adopt here. However this identification has a minor
effect on the symbol.

$\aPsDO^m(M; E, E)$
is closed under $*$ and so becomes a $*$-algebra.
we will use the abbreviation
$\aPsDO^m(E)=\aPsDO^m(M; E,E)$.

In particular,  for $m\le 0$,  $\aPsDO^m(E)$
is a \Fr $*$-subalgebra.

\subsubsection{Continuity under a Lie group action}

\begin{theorem}\label{thm:contcoordch}
Let $G$ be a Lie group acting on $M$.
For  $g\in G$ and $P\in \aPsDO^m(M; E,F)$,
let
\begin{equation}
g(P)=g_*(P)=g\circ P \circ g^{-1}
\end{equation}
where in the right hand side $g$ acts on
$\Gamma_c^\infty(M,E)$ and $\Gamma^\infty(M, E)$
by push-forward. Then this defines a
continuous action $G$ on $\aPsDO^m(M; E,F)$. 
\end{theorem}

\begin{proof}
By the topology we choose, it is sufficient to prove
the local version: without loss of generality,
we assume $P=a(x,D)$ is supported in an open
subset $X$ of $\R^n$. Let $\kappa_n$ be a sequence
of elements in $G$ converging to the identity,
we may assume all the $X_{\kappa_n}$ is contained
in a single open subset $X'$ of $\R^n$.

We need to show that under any semi-norm $p_{\alpha, \beta}$,
\[
p_{\alpha, \beta}(\kappa_n(P))\to p_{\alpha, \beta}(P)
\]

In the proof of the local invariance theorem \ref{locinv}
in \cite{Ho3} (theorem 18.1.17 there).
We choose $\phi\in C^\infty_c(X')$, such that $\phi(x)=1$
for any $x$ in the support of any $\kappa_n(P)$.

For any $\kappa=\kappa_n$
it is shown that
\begin{equation}
  a_\kappa(x,\eta)=\phi(x)e^{-i\langle \kappa(x), \eta\rangle}
  a(x,D)(\phi(x)
  e^{i\langle \kappa(x), \eta\rangle}
\end{equation}
is in $S^m$ and it is implicit in the proof that for all
$\alpha, \beta\in \N^n$, there is a constant $C_{\alpha, \beta}$
\begin{equation}
p_{\alpha, \beta}(a_{\kappa_n}(x,D)) \le C_{\alpha, \beta}
\end{equation}
which is uniform on $\{ \kappa_n \}$.

Each $a_{\kappa}(x,\eta)$ is shown to be an oscillatory
integral with the phase function 
\begin{equation}
  f_{x, \eta}(y,\xi)
  =\langle x-y, \xi\rangle -\langle \kappa(x)-\kappa(y), \eta\rangle
\end{equation}
with $\kappa$ converging to the identity. 
Before taking sup-norm over $x$ and $\eta$ to
get the semi-norms, the integrals
\begin{equation}
  (1+\abs{\eta})^{-m+\abs{\alpha}} D^\alpha_\eta D^\beta_x a_{\kappa}(x,\eta)
\end{equation}
are also absolutely bounded oscillatory integrals. The phase
functions and amplitude functions for these integrals
converge at every point $y,\eta$ as $\kappa\to id$.
Therefore these integrals are continuous with respect to
$\kappa$ in $G$.
\end{proof}

The continuity of $G$ action on the space $Op^m$ by
$g(P)$ is proved in \cite{AS1}.

\subsubsection{The distributional kernel}

\begin{proposition}\label{psdo:conti}
The map
\begin{equation*}
\aPsDO^m(M; E,F) \to \sL(\Gamma^\infty_c(M,E), \Gamma^\infty(M,F)),
\end{equation*}
from a symbol space to its action,
is continuous and extends to a bounded linear operator between
Sobolev spaces of the sections
\begin{equation}
P: H_{s}(M, E) \to H_{s-m}(M, F),
\end{equation}
for any real number $s$.
Moreover, the following natural map  is continuous:
\begin{equation}
\aPsDO^m(M;E,F) \to \sL_S(H_{s}(M, E), H_{s-m}(M, F))
\end{equation}
where in the right hand side ($\sL_S$ term)
is given the strong topology.
\end{proposition}

The above proposition follows from its local
version, see \cite{Ho3} for the proof.

We now recall the topology on $Op^m$ in \cite{AS1}.
\begin{definition}
Let $Op^m=Op^m(M;E,F)$ is the space of all continuous
linear maps in 
\[
\sL(\Gamma^\infty_c(M,E), \Gamma^\infty(M,F))
\]
which extend to 
\[
\sL_S(H_{s}(M, E), H_{s-m}(M, F))
\]
for all $s$. $Op^m$ is given the induced limit topology.
\end{definition}

$Op^m$ is a \Fr space. Operators in $Op^m$ do not necessarily
have the pseudo-local property, so they contain more than
pseudo-differential operators. By definition the natural
inclusion map
\begin{equation}
\aPsDO^m(M;E,F) \to Op^m(M; E, F))
\end{equation}
is continuous.

\begin{proposition}\label{eq:jointscgobal}
The composition map 
\begin{equation}
\aPsDO^{m_1}(M;E_1,E_2)\x \aPsDO^{m_2}(M;E_2,E_3)
\to \aPsDO^{m_1+m_2}(M; E_1, E_3)
\end{equation} 
is jointly continuous.
\end{proposition}

\begin{proof}
By the closed graph theorem for \Fr spaces
we need only to show separate
continuity. Without loss of generality, we assume
all three bundles $E_i$ are the same, say $E$.
Let $P_j\to 0\in \aPsDO^m(M; E,E)$ and
$P\in \aPsDO^{m'}(M; E, E)$,
for any $u\in C^\infty_c(M, E)$,
we have $P_j\circ Pu\to 0$ and $P\circ P_ju\to 0$,
which imply separate continuity.
\end{proof}

\subsubsection{Wave front sets of pseudo-differential operators}

In particular, a pseudo-differential operator
$P$ has a distributional kernel $K_P$.

\begin{theorem}
Let $P\in \aPsDO^m(M;E,F)$, then
\begin{equation}
WF'(P) \subset \{(x,x,\xi, \xi): x\in M, \xi\in T^*_xM, \xi\not=0\}.
\end{equation}
In other words, let
$K_P$ be the distributional kernel of $P$, then
\begin{equation}
WF(K_P) \subset \{(x,x,\xi, -\xi): x\in M, \xi\in T^*_xM, \xi\not=0\}.
\end{equation}
\end{theorem}

\begin{theorem}
For $P\in \aPsDO^m(M;E,F)$,
$u\in\Gamma^{-\infty}(E)$,
\label{psdoslocal}
\begin{equation}
WF(Pu) \subset WF(u)
\end{equation}
\end{theorem}

This property is called the micro-local, or
strong pseudo-local property.
Yet it can be improved if the
symbol of $P$ is smoothing in some directions.

\begin{definition}
Let $P: \Gamma^\infty_c(E)\to \Gamma^{\infty}(F)$
be a pseudo-differential operator.
The essential support $\Ess(P)$ of $P$
is the compliment in $T^*M\mz$
of the largest open conic subset of the
cotangent bundle on which the symbol
has order $-\infty$.
\end{definition}

In particular, $P$ is a smoothing operator
if and only if $\Ess(P)=\ep$.

\begin{proposition}
If $P\in \aPsDO^m(M;E,F)$ and $\Gamma$ is a closed
conic subset of $T^*M\mz$, the following are
equivalent:
\begin{itemize}
\item[1.] $P$ is of order $-\infty$ in $T*M\mz\sm \Gamma$;
\item[2.] $WF'(K_P) \subset \{(\xi,\xi), \xi\in \Gamma\}$;
\item[3.] for all $u\in \Gamma^{-\infty}(E)$,
  $WF(Pu) \subset \Gamma \cap WF(u)$.
\end{itemize}
\end{proposition}

So the conclusion in theorem \ref{psdoslocal}
can be improved to:
\begin{equation}
WF(Pu) \subset WF(u) \cap \Ess(P).
\end{equation}

\begin{definition}
When $E=F$, the characteristic set $char(P)$
of $P$, is the subset of $T^*M$ where
the principal symbol is not invertible
as bundle morphism.
\end{definition}


The following is a generalization of the regularity
theorem for the elliptic operators.

\begin{theorem}\label{thm:regularity}
(Regularity)
For $P\in \aPsDO^m(M;E,F)$,
$u\in\Gamma^{-\infty}(E)$,
\begin{equation}
WF(u) \subset WF(Pu) \cup char(P).
\end{equation}
\end{theorem}

In particular when $P$ is elliptic,
$char(P)=0_{T^*M}$
(but the essential support of $P$ is the
maximal, $T^*M\mz$);
so combining the previous two theorems we have
$WF(Pu)=WF(u)$, the elliptic
regularity property.

\subsubsection{Classical polyhomogeneous symbols}
  
A symbol $p\in S^m(\R^n)$ is called classical polyhomogeneous
if $p$ has an asymptotic expansion 
\begin{equation}
  p(x,\xi)\sim \sum_{j\in \N} p_j(x,\xi)
\end{equation}
where $p_j(x,\xi)$ is homogeneous of degree $m-j$ in $\xi$
for $\abs{\xi}]\ge 1$, that is, if $\abs{\xi}\ge 1$,
\begin{equation}
  p_j(x, \xi)={\abs{\xi}}^{m-j} p_j(x, \xi/\abs{\xi}).
\end{equation}
  
It is straightforward to extend the definition of classical
polyhomogeneous symbols to manifold and bundles, when we
may choose a Riemannian metric to define $\abs{\xi_x}$,
however, the set of classical polyhomogeneous symbols
does not depend on the choice of the metric.

The pseudo-differential operators with classical polyhomogeneous
symbols form a closed subalgebra of $\aPsDO^\infty$, which is also
closed under pull-back by diffeomorphism. 
This subalgebra is the one we are interested,
and we shall denote those classical homogeneous subspaces by
\begin{equation}
{\cPsDO}^\infty(M;E,F)=\aPsDO^\infty_{phg}(M;E,F).
\end{equation}

\section{Group action and pseudo-differential calculus}

\subsection{Transversal ellipticity}

\begin{definition}
\label{def:fun}
For $X\in \g$,  let $X_M$ denote the
{\it fundamental vector field} generated by
the action of the one parameter group
corresponding to $X$:
\begin{equation}
X_M f(x)=
\left.\frac{d}{dt}\right\vert_{t=0}
f\left(e^{-tX}x\right)
\end{equation}
for any smooth function $f$ on $M$ and $x\in M$.
\end{definition}

The $G$-action decides a map $X\mapsto X_M$
from $\g$ to $\Gamma(TM)$.

Similarly, for a $G$-equivariant bundle $E$ over $M$,
the bundle version is a first order differential
operator determined by the flow of action generated
by $X$:
\begin{equation}
(X_E s)(x)=
\left.\frac{d}{dt}\right\vert_{t=0}
\left(\rho(e^{-tX})s\right)
\left(e^{-tX}x\right).
\end{equation}
Note that the $*$-operation on $\cPsDO^m(E,F)$
is preserved when the action of $G$ preserves the volume form.

\begin{definition}
Let
\begin{equation}
\mu_*^{M,G}: T^*M \to \g^*
\end{equation}
be the adjoint of the group action map
defined as: for any $(x, \xi)\in T^*M$
\begin{equation}
\pair{\mu_*^{M,G}(x, \xi), X}
=\pair{\xi, X_M(x)}.
\end{equation}
\end{definition}

The kernel for the moment map is usually denoted by
\begin{equation}
T^*_G M=(\mu_*^{M,G})^{-1}(0)=
\{(x, \xi) \in T^*M |
\pair{\xi, X_M} = 0 \;\forall X \in \g\}.
\end{equation}
$T^*_G M$ is a $G$-invariant conic closed
subspace of $T^*M \mz$.

\begin{definition}
Let $P\in \cPsDO^k(E,F)$ be
a $G$-invariant pseudo-differential
operator. $P$ is called
{\it transversally elliptic}
relative to the action of $G$, if its principal symbol,
well defined as an element
\begin{equation}
\sigma_{P}: S^m(T^*M; Hom(E,F))/S^{m-1}(T^*M; Hom(E,F))
\end{equation}
is invertible for $(x, \xi) \in T^*M\mz$ as an element in
$Hom(E,F)$.

$P$ is called strongly transversally elliptic
if there exists a conic neighborhood $U$
of $T^*_G M\mz$ in $T^*M\mz$ and
an inverse principal symbol
\begin{equation}
\sigma_Q\in S^{-m}(U; Hom(F,E))/S^{-m-1}(U; Hom(F,E)).
\end{equation}
\end{definition}

The strong transversal ellipticity condition
is automatic for transversally elliptic classical
polyhomogeneous symbols, which is our primary concern.

\subsection{Crossed product algebra and action groupoid}

\begin{definition}\label{def:prod}
$\A=C^\infty_c(M) \rx G$ is defined
to be the $*$-algebra of
smooth crossed product. That is,
$\phi\in \A$
means $\phi\in C^\infty_c(M\x G)$
with product and adjoint:
\begin{equation}
\phi*\psi (x, g)
=\int_G
\phi(x, h)
\psi(h^{-1}x, h^{-1}g)
d\mu(h)
\end{equation}
\begin{equation}
\phi^*(x,g)=\overline{\phi(g^{-1}x, g^{-1})}.
\end{equation}
\end{definition}

To show that the above product forms a $*$-algebra,
the only nontrivial part is the associativity of
the product, which can be computed directly.
Instead we will only show that
this binary operation 
is a special case
of the convolution product of
the action groupoid. This fact also
supports our similar arguments later on.

A concise definition of a groupoid $\G$
is that a groupoid is a small category
(in which all morphism form a set)
in which every morphism is invertible. 
A more detailed definition (such as in
\cite{MRW:groupoidcstar}) is
\begin{equation}
\G=\left(\G_1, \G_0; \tau, \sigma, \iota, \cdot, ()^{-1}\right)
\end{equation}
with the five maps satisfying the axioms for
a small category. We elaborate only on our main example.

\begin{example}
The action groupoid is a groupoid $\G_\rho$,
where the set of unities $\G_0$ (or objects) is $M$,
and the set of arrows $\G_1$ (or morphism) is $M\x G$,
and with the five maps as follows:

(1) The target map $\tau: M\x G \to M$: $\tau(x,g)=x$.

(2) The source map $\sigma: M\x G \to M$: 
$\sigma(x,g)=\rho(g^{-1})x=g^{-1}x$
(we will sometimes omit the action $\rho$ if no
confusion will be caused).

(3) The unity map $\iota: M \to M\x G$: $\iota(x)=(x, e)$.

(4) The partially defined multiplication defined on the subset
of $\G_1\x \G_1$ where the source of the first component matches
the target of the second component, is given simply by the 
multiplication of the group $G$.  That is,
$(y, h)\cdot (x, g)$ is defined only when $y=hx$
and
\begin{equation}
(hx, h)\cdot (x,g) = (hx, hg). 
\end{equation}

(5) The inverse map 
$()^{-1}: M\x G \to M\x G$: $(x,g)^{-1}=(g^{-1}x, g^{-1})$.

\end{example}

Remark: our definition of the action groupoid
is slightly different from another commonly used in the literature,
where typically an arrow $(x,g)$ has source $x$ and target
$gx$. The map $(x,g)\mapsto (g^{-1}x, g)$
gives a groupoid isomorphism between them.

\begin{definition}
A groupoid
$\G=\left(\G_1, \G_0; \tau, \sigma, \iota, \cdot, ()^{-1}\right)$
is a topological groupoid if all the following are true 
(1) $\G_1$, $\G_0$ are topological spaces;
(2) all five maps are continuous;
(3) $\tau$, $\sigma$ are open maps;
(4) $\iota$ is a homeomorphism onto its image.
A topological groupoid is called
locally compact if $\G_1$ is locally compact
(so is $\G_0$ as a consequence).
$\G$ is called a smooth groupoid when all of the following are
true: (1) $\G_1$, $\G_0$ are both smooth manifolds; (2)
all five maps are smooth; (3) $\tau$, $\sigma$ are both submersions;
(4) $\iota$ is an embedding.
\end{definition}

The action groupoid $\G$ is smooth.

For a locally compact topological groupoid, we recall
the definition of a Haar measure
which gives a convolution algebra.

\begin{definition}
A smooth Haar system on a smooth groupoid $\G$
is a family of positive Radon measures
$\lambda^\bullet=\{\lambda^x: x\in \G_0\}$
on $\G_1$ satisfying the following conditions:

(1) For any $x\in \G_0$, the support of $\lambda^x$ is
in $\G^x=\{\alpha\in \G_1: \tau(\alpha)=x\}$ and
$\lambda^x$ is a smooth measure on $G^x$.

(2) (Left invariance) For any $x\in \G_0$
and any continuous function
$f: \G^x\to \C$
and any $\alpha\in \G^x$ we have
\begin{equation}
\int_{\G^x} f(\beta) d\lambda^x(\beta)
=\int_{\G^{\sigma(\alpha)}} 
f(\alpha\beta)d\lambda^{\sigma(\alpha)}(\beta)
\end{equation}
(in other words, 
$\alpha_*(\lambda^{\sigma(\alpha)})=\lambda^{\tau(\alpha)}$).

(3) (Smoothness) For any $\phi\in C^\infty_c(\G)$,
the map 
\begin{equation}
x\mapsto \int_{\G^x} \phi(\beta) d\lambda^x(\beta)
\end{equation}
is a smooth function on $\G_0$.
\end{definition}

\begin{proposition}\label{prop:conv}
A smooth Haar system on a smooth groupoid $\G$
defines a convolution product on $\C^\infty_c(\G)$ by
\begin{equation}
(\phi*\psi)(\alpha)=\int_{\G^{\tau(\alpha)}}
\phi(\beta) \psi(\beta^{-1}\alpha) d\lambda^{\tau(\alpha)}.
\end{equation}
And it is a $*$-algebra with
\begin{equation}
\phi^*(\alpha)=\overline{\phi(\alpha^{-1})}.
\end{equation}
\end{proposition}

For an action groupoid $\G_\rho$,
let $\lambda$ be the Haar measure on $G$, 
$\G^x=\{(x,g): g\in G\}$ is homeomorphic
to $G$ by the projection $\pi$ on the second
component. So we simply
define $\lambda^x=\pi^*\lambda$.
In other words, $x$ is nothing but a label for the
copy of $G$. Left invariance of $\{\lambda^x: x\in M\}$
simply follows from that of $G$. The rest of the
properties in the
Haar system definition are also easy to check.

In proposition \ref{prop:conv}
let $\alpha=(x,g)$, $\beta=(x,h)$
(since we need $\beta\in \G^{\tau(\alpha)}$),
we have $\beta^{-1}\alpha=(h^{-1}x, h^{-1}g)$.
It is clear that definition \ref{def:prod}
of the product
is just a special case of definition \ref{prop:conv}
of the convolution product.

\begin{lemma}\label{lem:Amod}
Let $E$ be a $G$-bundle over $M$.
Denote by $\rho (g): E_x \to E_{gx}$ the
action of an element $g\in G$.  Then $\Gamma(M, E)$ is a
$\A$-module, with
the actions $\phi\in \A$
denoted by $\rho$, for any $s \in \Gamma(M, E)$,
\begin{equation}
(\rho (\phi)s)(x) = \int_G \phi(x, g)
\rho (g)(s(g^{-1} x)) d\mu(g).
\end{equation}
\end{lemma}

For the proof, instead of direct verification,
we relate it to some standard results
of groupoid action in the following discussion.

\begin{definition}
Let $\G$ be a smooth groupoid,
a left $\G$ action
a smooth manifold $M$ is a pair
$(\rho, t)$, where $t: M\to \G_0$
is a smooth map and
\begin{equation}
\rho: \G_1 \x_{\G_0} M \to M
\end{equation}
where
\begin{equation}
\G_1 \x_{\G_0}M=\{(\alpha, x)\in \G_1\x M: \sigma(\alpha)=t(x)\},
\end{equation}
with the following properties

(1) for any $z\in \G_0$, $\iota(z)=t(x)$, we have $\rho(\iota(z),x)=x$,

(2) for all $\alpha$, $\beta$ with $\sigma(\alpha)=\tau(\beta)$
and $\sigma(\beta)=t(x)$, we have
\begin{equation}
\rho(\alpha\beta, x)=\rho(\alpha, \rho(\beta, x)).
\end{equation}

When $\rho$ is clear under context, we 
often use the abbreviations $\rho(\alpha)x$,
or even $\alpha x$ for $\rho(\alpha,x)$.

\end{definition}

Let $E$ be a hermitian bundle over $M$ and
let $\rho$ be a unitary $G$ action on $E$,
which means for each $x\in M$,
\begin{equation}
  \rho(g): E_x \to E_{gx}
\end{equation}
is a unitary. 

$\rho$ introduces a $\G_{\rho}$
action on $U(E)$ (the bundle of fiber-wise
unitary transformations on $E$) in the
following way.
The groupoid action is given the pair
$(\tilde\rho, t)$ where
$t: U(E) \to M$ is the projection along the
fiber, and $\tilde\rho$ is defined as
\begin{equation}
\tilde\rho((gx,g), u_x) = \rho(g) u_x \in U(E_{gx})
\end{equation}
for $u_x\in U(E_x)$.

\begin{proposition}
Suppose $\G$ acts on $M$ from the left,
$C^\infty_c(M)$ is a left
$C^\infty_c(\G)$-module by the following action:
for $\phi\in C^\infty_c(\G)$ and $f\in C^\infty_c(M)$,
\begin{equation}
(\phi f)(x)=
\int_{\G^x} \phi(\alpha) f(\alpha^{-1} x) d\lambda^x(\alpha).
\end{equation}
\end{proposition}
See \cite{MRW:groupoidcstar} for a proof.     

From the natural action of $C^\infty_c(U(E))$
on $\Gamma(M,E)$, and combined with the action
of $\A$ on $C^\infty_c(U(E))$, 
the action in lemma \ref{lem:Amod}
is shown to be the composition.

\begin{lemma}\label{bdd2}
Let $E$ be a $G$-vector bundle, for any
$\phi\in \A$, the linear continuous operator 
\begin{equation}
\rho (\phi) : \Gamma^\infty(E) \to \Gamma^\infty(E)
\end{equation}
extends to a bounded linear operator on any Sobolev space,
that is, for any $s\in \R$,
\begin{equation}
\rho (\phi) : H_s(E) \to H_s(E)
\end{equation}
is bounded.
\end{lemma}
The proof of the lemma is straightforward.

\subsection{Crossed product with pseudo-differential operators}

For a manifold with a group action,
we will often use operators which are not
in the algebra of pseudo-differential
operators. 
For example, in our purpose we frequently use
the commutator of a pseudo-differential operator and
the action of $\rho(\phi)$, which is in general not
pseudo-local.
We will define a large enough algebra
$\cPsDO^k(E,E)\rx G$ which is, in short,
the one generated by both the group
action and pseudo-differential operators.
$\cPsDO^k(E,E)\rx G$ may also be viewed
as groupoid algebra
as introduced in the previous section.
Here the main advantage of the groupoid
viewpoint is that many properties
of the new algebra can be easily reproduced.

Similar to definition \ref{def:prod}, we introduce:
\begin{definition}\label{xpsdo}
Let $E$ be an hermitian vector bundle over $M$. We
define $\cPsDO^\infty(E,E)\rx G$ as the algebra of families
of pseudo-differential operators $P(g)\in \cPsDO^k(E,E)$,
and the product of $P(g)$ with $Q(g)\in \cPsDO^l(E,E)$
is a family of pseudo-differential operators in 
$\cPsDO^{k+l}(E,E)$,
\begin{equation}
(P*Q)(g)=\int_G P(h) \cdot
\left[ ((h^{-1})_*Q)(h^{-1}g)\right]
d\mu(h),
\end{equation}
\begin{equation}
P^*(g)=\overline{g^{-1}P(g^{-1})^*g}.
\end{equation}
\end{definition}

Since the composition, 
the conjugation with $g,h\in G$, and 
the integration on the parameter space 
are closed on $\cPsDO^{k+l}(E,E)$, and the continuity
of the integrand over $G$ ensure that the
integral (as a Riemann integral over a Lie
group with value in a \Fr space) is still
in $\cPsDO^k(E,E)\rx G$. We may show
$\cPsDO^\infty(E,E)$ and $\cPsDO^0(E,E)$ are algebras
by the groupoid argument over a \Fr algebra.
 
It follows directly that $\A\subset \cPsDO^0(E,E)\rx G$;
$\cPsDO^0(E,E)\rx G$ is not a unital algebra;
and $\cPsDO^\infty(E,E)\rx G$ 
is a $\cPsDO^\infty(E,E)$-bimodule.

In fact, we are only concerned with the following
representation of the action of $\cPsDO^k(E,E)\rx G$
on smooth sections of $E$, and 
$\sH=L^2(E)$ or in general, any Sobolev space $H^s(E)$.
We use the same notions as in lemma \ref{lem:Amod}.

\begin{proposition} There is a natural continuous action of 
$\cPsDO^0(E,E)\rx G$ on $\Gamma^\infty(E)$ defined
as follows. For $P=P(g)$ in $\cPsDO^0(E,E)\rx G$,
\begin{equation}
  (Ps)(x)=
  \int_G P(g) \rho (g)(s(g^{-1} x)) d\mu(g),
\end{equation}
The action is extendable to $H_s(E)$ for any
$s\in\R$, so the action gives an element in $Op^0$.

In fact, for any $s\in \R$,
an element of $\cPsDO^k(E,E)\rx G$ extends to
$H^s(E)$ to $H^{s+k}(E)$, defining an element in $Op^k$.
\end{proposition}

\begin{proof}
By Lemma \ref{bdd2} and the uniform continuity
of the family $P=P(g)$ over $g\in G$, and 
the continuity of integration on $G$,
we may routinely check all the above conclusions,
based on the corresponding properties of
pseudo-differential operators. 
\end{proof}

From this point on we use the same notation for an element of
$\cPsDO^k(E,E)\rx G$ and its representation in $Op^k$.

An ordinary pseudo-differential operator $Q$,
viewed as a constant family over $G$.
But the algebra of pseudo-differential
operators (as constant families) is not
a subalgebra of $\cPsDO^\infty(E,E)\rx G$:
the composition of two such operators is given
by the *-product, --- just as convolution of functions
is different from pointwise multiplication.

The following compositions between
pseudo-differential operators and elements
in $\cPsDO^\infty(E,E)\rx G$ are induced
from their actions on $\Gamma^\infty(E)$:
for $P=P(g)\in\cPsDO^k(E,E)\rx G$
and $Q\in\cPsDO^k(E,E)$
\begin{equation}
\begin{split}
(P*Q)(g)&=P(g)\cdot g^*(Q)\\
(Q*P)(g)&=Q \cdot P(g).
\end{split}
\end{equation}

It is straightforward to check that
$\cPsDO^\infty(E,E)\rx G$ is a graded bimodule
of pseudo-differential algebra $\cPsDO^\infty(E,E)$:
\begin{equation}
\begin{split}
(\cPsDO^k(E,E)\rx G)\cdot\cPsDO^l(E,E)&\subset \cPsDO^{k+l}(E,E)\rx G\\
\cPsDO^l(E,E)\cdot(\cPsDO^k(E,E)\rx G)&\subset \cPsDO^{k+l}(E,E)\rx G.
\end{split}
\end{equation}

\subsection{Wave front sets of crossed product algebras}

\begin{lemma}\label{lemma:xpsdo}
For any $P\in \cPsDO^k(E,E)\rx G$,
\begin{equation}
WF'(P) \subset
\left\{
(g_*\xi_x, \xi_x): \,
(x,\xi) \in \Ess(P(g)),\; \xi_x\in (T_G^*M)_x
\right\}
\end{equation}
\end{lemma}

\begin{proof}
Recall that under the notation,
$\rho: M\x G \to M$ is the action of $G$ on $M$,
$\rho^*E$ is the pull-back bundle of $E$ on $M\x G$.
As an operator, $P$ is the composition
of
\begin{equation}As(x, g)=  \rho (g)(u(g^{-1} x)),
\end{equation}
the family $P(g)$, and
\begin{equation}
Bu(x)=\int_G u(g, x) d\mu(g).
\end{equation}
We now show that
$A=(f_A)_*$ and $B=(f_B)_*$ are both push-forward operators,
an embedding and a submersion respectively.
$f_A: M\x G \to M\x M\x G$
\begin{equation}
f_A(x,g)=(gx, x,g)
\end{equation}
which is the lower part of the bundle map
\begin{equation}
\tilde{f}_A: \rho^* E \to E\boxtimes\rho^*E
\end{equation}
\begin{equation}
\tilde{f}_A: e_{(x,g)} \mapsto \rho(g)e_{gx};
\end{equation}
$f_B: \rho^*E \to E$
is the projection to the first
component, which is the lower part of
the bundle map
\begin{equation}
\tilde{f}_B: E\boxtimes\rho^*E \to E
\end{equation}
\begin{equation}
\tilde{f}_B: (e_y, e'_{x,g}) \mapsto e_y.
\end{equation}

It is easy to check that $A=(f_A)_*$
and $B=(f_B)_*$.

By theorem \ref{thm:push2}
we have
\begin{equation}
WF'(A)\subset\{(g_*\xi_x, \xi_x,\gamma_g):
\mu_*(\xi_x)=\gamma_g,
\}
\end{equation}
where $\mu_*=\mu^{M,G}_*$
is the moment map (with
kernel $\{\xi_x\in (T_G^*M)_x\}$)
defined earlier.
In particular both $WF'_M(A)$ and
$WF'_{M\x G} (A)$ are empty.
And for the wave relation of
$P=P(g)$ we apply theorem \ref{psdoslocal},
which has empty projections.
At last by theorem \ref{thm:push2} we have
\begin{equation}
WF'(B)\subset \sN^*_{M\x M}(M\x M\x G)
=\{(\eta_y,\xi_x,0)\}.
\end{equation}
Although it has nonempty right projections,
the intersection with the left projection
of $P$ is empty.
Apply theorem \ref{thm:compo}
twice we will reach the conclusion.
\end{proof}

\begin{corollary}\label{lemma:xwf}
For any $\phi \in\A$,
\begin{equation*}
WF'(\rho(\phi)) \subset
\left\{
(g_*\xi_x, \xi_x): \,
(x,g) \in supp({\phi}),\; \xi_x\in (T_G^*M)_x
\right\} .
\end{equation*}
\end{corollary}

\subsection{Transversally smoothing operators}

\begin{definition}
Let $\sK_G \subset \cPsDO^\infty(E,E)$ be
the set of pseudo-differential operators
that annihilate $\cPsDO^\infty(E,E)\rx G$
modulo $\sK$, the ideal of smoothing operators,
with composition product. We shall call
elements of $\sK_G$ transversally smoothing
operators.
\end{definition}

$\sK_G$ can be described
in terms of conditions on symbols.
It contains those pseudo-differential operators
with asymptotically zero symbol on $T^*_GM$.
In other words, for $P\in \cPsDO^\infty$, if
\begin{equation}
 \Ess(P)\cap T^*_G(M)=\ep
\end{equation}
then $P\in \K_G$.

So $\sK_G$ is an ideal if $\cPsDO^\infty$,
containing $\sK$ and some pseudo-differential
operators of arbitrarily high order.

\begin{example}\label{eg:101} 
Fix a pair of conic neighborhood $U_1$ and
$U_2$ of $T^*_G$, such that $U_1/\R_+$
is relatively compact in $U_2/\R_+$.
There exists a real valued, positive,
smooth function $\chi$ on $T^*M\mz$ satisfying
the following properties: (1) $\chi(\xi)=\chi(\xi/\abs{\xi}))$
for $\abs{\xi}\ge 1$, that is, $\chi(\xi)$ as a symbol
is homogeneous of degree $0$ for $\abs{\xi}\ge 1$;
(2) $\chi(\xi)= 0$ if $\xi\in U_1$ and $\abs{\xi}\ge 1$;
(3) $\chi(\xi)= 1$ if $\xi\not\in U_2$;
(4) $\chi(\xi)= 1$ if $\abs{\xi}\le 1/2$;
(5) $\chi$ is $G$ invariant.

Let $P_\chi$ be a scalar pseudo-differential
operator with symbol $\chi$. Since 
$\chi=\sigma(P_\chi)$ vanishes on a conic
neighborhood of $T^*_GM$,
$P_\chi\in\sK_G$.
\end{example}

More transversally smoothing operators can be
constructed through this example.

\begin{example}\label{eg:102} 
Let $P\in \cPsDO^m(M;E,F)$ be any pseudo-differential
operator of order $m$. Then by the symbol expansion
formula $P_\chi P\in \sK_g$ and it generically has order $m$.
Same is true for $P P_\chi$ and $P_\chi P P_\chi$.
\end{example}

\subsection{Transversal parametrix}

\begin{proposition}
\label{prop:parametrix} 
Let $P\in \cPsDO^m(M;E,F)$ be transversally
elliptic. Then there exists a $Q\in \cPsDO^{-m}(M;F,E)$
with principal symbol
$\sigma(Q)=\sigma(P)^{-1}$ on $T^*_G M$,
such that $K=1_F-PQ$ and $K'=1_E-QP$
are both transversally smoothing.
\end{proposition}

For brevity we will call such a $Q$
a {\it transversal parametrix} for $P$.

\begin{proof}
The construction of the transversal
parametrix repeats essentially
the construction of the parametrix
of an elliptic pseudo-differential
operator( see \cite{BEP}).

In the scope of this proof,
we use the sign ``$\equiv_V$'' for the
asymptotic equivalence (up to a difference of
degree $-\infty$) on a conic neighborhood $V$
of $T^*_GM\mz$. We omit $V$ when it is $T^*M\mz$.

Let $\sigma(P)$ be the symbol of $P$,
with order $m$ and 
let $Q_0$ be a pseudo-differential operator
of order $-m$ and with inverse symbol
on a conic neighborhood $V$ of $T_G^*M$.
Then $R_0=Q_0P-1$  is a pseudo-differential operator,
in general not necessarily of negative order,
but it has negative ordered symbol
on a conic neighborhood of $T^*_GM$.
By composition with an operator of
the form $P_{1-\chi}$ as in example
\ref{eg:101} with $U_1\subset$ $U_2$
$\subset$ $V$,
we get a symbol for a pseudo-differential operator
$R$ of negative order, and
\begin{equation}
\sigma(R)\equiv_{U_1}\sigma(R_0).
\end{equation}
Let $C$ be a pseudo-differential operator
with the following asymptotic expansion
of its symbol: 
\begin{equation}
\sigma(C) \sim \sum\limits_{i=0}^{\infty}(-1)^i\sigma(R^i)
\end{equation}
where $\sigma(R^i)$ is the symbol of
the $i$-$th$ power $R^i$ of $R$.
Thus
\begin{equation}
\sigma(C(1+R))\equiv \sigma(1).
\end{equation}

We claim that $Q=CQ_0$ is a pseudo-differential
with properties we need. First we check the
properties about $1-QP$.
We have
\begin{equation}
\sigma(QP)\equiv \sigma(C(1+R_0))
\equiv\sigma(C(1+R))
\equiv_{U_1} \sigma(1).
\end{equation}

To prove the properties about $1_F-PQ$,
we may construct $Q'$ a similarly so
$1-PQ'\equiv_{U_1} 0$. But with the
existence of such a $Q'$ we know $Q$
work in place of $Q'$:
\begin{equation}
\sigma(Q)-\sigma(Q')\equiv_{U_1}
\sigma(QPQ')-\sigma(QPQ')=0,
\end{equation}
therefore
\begin{equation}
\sigma(1-PQ)\equiv_{U_1}
\sigma(1-PQ')+\sigma(P(Q'-Q))
\equiv_{U_1} 0.
\end{equation}

\end{proof}

Consequently, $1-PQ$ and $1-QP$
composed with any operator in $\cPsDO^\infty_G$
is smoothing. In particular, their composition
with any $\rho(\phi)$ for $\phi\in \A$ is
smoothing.

\subsection{Summary of algebraic properties}

We summarize some facts about the
algebra $\cPsDO^k(E,E)\rx G$ and the
pseudo-differential algebra $\cPsDO^k(E,E)$.

\begin{enumerate}
\item $\cPsDO_G^k\subset \cPsDO_G^l$ for $k\le l$;

\item $\cPsDO_G^k\cdot \cPsDO_G^l\subset \cPsDO_G^{k+l}$;

\item $\cPsDO^k\cdot \cPsDO_G^l\subset \cPsDO_G^{k+l}$; 
\quad $\cPsDO_G^k\cdot \cPsDO^l\subset \cPsDO_G^{k+l}$;

\item $\cPsDO_G^0\subset B(\sH)$;

\item When $r>0$, $A\in \cPsDO_G^{-r}$ is compact;

\item $\cPsDO_G^k\cdot \sK_G\subset \sK$; 
\quad $\sK_G \cdot \cPsDO_G^k\subset \sK$;

\item $\cPsDO_G^{-\infty}\subset \sK$.


\end{enumerate}

\section{Trace formula and spectral analysis}

\subsection{Trace formula on $\cPsDO_G$}

For $\sH=L^2(E)$, we recall that $\sL^1(\sH)$ is
the trace class operators and for $A\in\sL^1(\sH)$,
$\norm{A}_1$ is the trace class trace class norm.

\begin{proposition}\label{prop:tr:k}
If $k<-\dim{M}-1$, then the extension from
smooth sections to $L^2$ sections sends
\begin{equation}
  \cPsDO^k(E,E)\rx G \subset \sL^1(\sH).
\end{equation}
Moreover the above embedding 
from the the \Fr algebra $\cPsDO^k(E,E)\rx G$
to the Banach space $(\sL^1(\sH), \norm{\cdot}_1)$
is continuous.
\end{proposition}

\begin{proof}
We take a self-adjoint, elliptic, invertible pseudo-differential
operator, such as the second order invariant Laplacian
$1+\Delta$ on scalar functions $C^\infty(M)$.
$A(1+\Delta)^{-k/2}$ is bounded,
and $(1+\Delta)^{k}$ is trace class by the Weyl's formula:
\begin{equation}
\mu_n(1+\Delta) = c_n n^{\dim{M}/2} (1+o(1/n))
\end{equation}
where the positive $c_n$ depends only on $n$ and the volume of $M$.

To show the trace class norm is bounded on $\Psi_G^{k}$,
we can use
\begin{equation}
\norm{A}_1 \le \norm{A(1+\Delta)^{-k/2}}_{B(\sH)} 
\norm{(1+\Delta)^{k/2}}_1,
\end{equation}
together with joint continuity of the multiplicative
map in $\cPsDO$.
\end{proof}

The above argument gives a stronger result than
the following one by just looking at the distributional kernel.
For $k$ sufficiently negative, the distributional
kernel of $A\in \cPsDO^k_G$ is $C^r$ for all
$r< -k-\dim{M}$. Stinespring
\cite{Stinespring} showed that if an integral operator
is trace class if its kernel is $C^r$ with
$r \ge [\dim{M}/2]+1$.
However in proposition \ref{prop:tr:k} the condition
on $k$ can be relaxed in most interesting cases to
\begin{equation}
  \label{eq:tr:mg}
  k < -\dim{M_0/G}
\end{equation}
where $M_0$ is the union of all of principal
orbits of the $G$-action, an open dense subset of $M$.
(see also \cite{BrHe}).

Also for $P\in \cPsDO^{k}_G$, $k<-\dim{M}-1$,
we have a well defined Fredholm determinant
\begin{equation}
  \det{(I+\lambda P)} \in \C
\end{equation}
for any $\lambda\in \C$ and is an entire function
on $\lambda\in\C$. The Fredholm determinant
is a useful tool in the study of trace and
eigenvalue related problems. We refer
to \cite{GGK} for details. However we do not use
the Fredholm determinant here, except for emphasizing
its role in the proof of the Lidskii's trace theorem.

\begin{proposition}
Assume  $k<-\dim{M}-1$.
For $P=\{P(g)\}\in \cPsDO^k_G$ let
\begin{equation}
K_P(x,y) \in \Gamma^\infty(E\boxtimes E)
\end{equation}
be the continuous representation of the
distributional kernel of $P(g)$.
Then
\begin{equation}\label{eq:tr:inv}
Trace(P)=\int_G\int_M tr_x(K_P(x, gx)\rho_x(g))dvol(x)d\mu(g),
\end{equation}
where $tr_x$ is the fiber-wise trace on $E_x\otimes E_x$.
\end{proposition}

\begin{proof}
By the continuity of the trace under the trace class norm,
we easily reduce the statement to its version on $\cPsDO^k$:
for $R\in \cPsDO^k(E,E)$ with distributional kernel
\[
K_R(x,y) \in \Gamma^\infty(E\boxtimes E)
\]
which is continuous,
\begin{equation}\label{eq:318}
Trace(R)=\int_M tr_x(K_P(x,x))dvol(x).
\end{equation}

To prove (\ref{eq:318}), we first prove it for
$R\in\cPsDO^{-\infty}$. There (\ref{eq:318})
holds since the set of operators of finite rank
is dense under the topology defined by the trace class norm.
Next we observe the right hand side of (\ref{eq:318})
is a unique continuous extension of the trace from
$\cPsDO^{-\infty}$ to $\cPsDO^k$.

By the Lidskii trace theorem \cite{GGK}, any
such extension of the trace on operators of
finite rank must be equal to
\begin{equation}
Trace(R)=\sum{\mu_n(P)}
\end{equation}
where $\mu_n(P)$ is the $n$-$th$ non-zero
eigenvalue of $P$ counting multiplicity
(for a compact operator the nonzero eigenvalues
have finite multiplicity); so it is the unique
trace. The harder part in the proof of the Lidskii's
theorem is to show that on operators with no eigenvalue,
any reasonably continuous extension of the trace vanishes.

\end{proof}

\subsection{Asymptotic analysis of oscillatory integrals}
\label{sec:32}

\subsubsection{Localization of trace formula}

In the trace formula (\ref{eq:tr:inv}) it
is not convenient to relate the distributional
kernels to the symbols of the operators.
We try to decompose the global formula 
into a sum of terms which are all local
in coordinate charts.

To do this we first choose a finite atlas of $M$,
on each chart in the atlas the bundles are trivial.
We take a refinement of this atlas so that
each new chart has a domain which is
relative compact subset of the domain of any
chart in the old atlas.
Then we use a partition of
unity on $G$ fine each so that we can assume
on each small piece the $G$ action on any
of the domains of the new atlas.
Next by a partition of unity on $M$, we
can decompose (\ref{eq:tr:inv}) into the sum of
integrals on open subset $U\x G$ such that
$g(U)$ is contained in the same chart for
all $g$ in the support of $g\mapsto P(g)$.

Locally the kernel of a pseudo-differential operator
$Q(g)$ is of the form
\begin{equation}
  K_{Q(g)}(x,y)=(2\pi)^{-n}\int_{\R^n}
  e^{i(x-y)\cdot \xi} q(x, \xi, g) d\xi,
\end{equation}
in particular,
\begin{equation}
  K_{Q(g)}(x,x)=(2\pi)^{-n}\int_{\R^n}
  q(x, \xi, g) d\xi.
\end{equation}

But for $Q=\{Q(g)\}$, the kernel on the
diagonal is of the form:
\begin{equation}\label{eq:os322}
  Q(x,x)=\int_G\int_{\R^n}
  e^{i(x-gx)\cdot \xi} p(x, \xi, g)  d\xi dvol(x) d\mu(g)
\end{equation}
and it is an oscillatory integral with phase function
$(x-gx)\cdot \xi$.

In the case there is no group action, equivalently
$G=\{e\}$, we find that the phase function
\begin{equation}
\phi(x,g,\xi)= \pair{x-gx, \xi}
\end{equation}
in formula (\ref{eq:os322}) is zero hence no
oscillatory integrals is directly
involved in the kernel.

\subsubsection{Asymptotic expansion of trace}

Here we use a simple example to show the need to
study the asymptotic behavior of oscillatory
integrals. 

We take the scalar functions on $\R^n$. 
Given a strictly positive, elliptic
differential operator $P$, so that the symbol
$p(x,\xi)$ is homogeneous of degree $m>0$.
We shall need to analyze the behavior of the
resolvant $(P-\lambda)^{-1}$ for $\lambda$
on the negative half of the real axis,
and let $\abs{\lambda}$ tends to $\infty$.

If $m>n$, the resolvant is trace class for
any $\lambda\in (-\infty,0)$, since it is
a pseudo-differential operator of order $-m$.

\begin{equation}
  \label{eq:ex:100}
  Trace((P-\lambda)^{-1})
  =(2\pi)^{-n}
  \int_{\R^{2n}} e^{i(x-gx,\xi)} (p(x,\xi)-\lambda)^{-1} \;d\xi\, dx.
\end{equation}

Suppose we want to study
how the trace varies with respect to $\lambda$.
For convenience, 
let $\lambda=-\mu^m$ where $\mu$ is positive,
let $t$ be $\abs{\xi}$,
and let $\xi_1=\xi/t$ be the unit vector.

As in \cite{Seeley} and \cite{Shubin}, we treat the symbol
\begin{equation}
  \label{eq:ex:101}
  (p(x,\xi)-\lambda)^{-1}=  (t^m p(x,\xi_1)+\mu^m)^{-1}
\end{equation}
as a homogeneous symbol in $(\xi, \mu)$. With a substitution
$t=\mu s$ and then rename $\xi=s\xi_1$ the trace integral looks like
\begin{equation}
  \label{eq:ex:102}
  (2\pi)^{-n} \mu^{-m}
  \int_{\R^{2n}}
  e^{i\mu(x-gx,\xi)} (p(x,\xi)+1)^{-1} \;d\xi\, dx.
\end{equation}

Note that
$\mu=(-\lambda)^{1/m}$ appears only outside of the integral and
as a linear factor of the phase function. This showcase
demonstrates why the asymptotic behavior of oscillatory
integral appears.
It also shows when $G$ action is trivial there is
no oscillatory integral involved.

\subsection{Asymptotic expansion of oscillatory integrals}

\subsubsection{The question}

We are interested in integrals of the 
form (\ref{eq:os322}) with $p(x, g, \xi)$ polyhomogeneous
in $\xi$ and with an extra parameter $\mu$.
The asymptotic behavior when $\abs{\xi}$ and $\mu$ are large
is a major concern, so we formulate the abstract
oscillatory integral as
\begin{equation}\label{eq:os324}
I_{\phi,u}(\tau)=\int_{\R^N} e^{i\tau \phi(y)} u(y) dy
\end{equation}
where $\tau>0$, $u$ is smooth.
In practice, $\phi$ and $u$ may also contains parameters.
The question is to decide the asymptotic behavior of
the integral (\ref{eq:os324}) as $\tau$ tends to $\infty$.

The derivative of the phase function $\phi'$
(with respect to $y$) matters the most in the
pattern of the possible asymptotic expansion.

The first basic fact is that $\phi'$ is away
from zero, then in (\ref{eq:os324})
\begin{equation}
  \label{eq:os:decay}
  I_{\phi, u}(\tau)=o(\tau^{-N})
\end{equation}
for any $N>0$.

The next basic fact is that if $\phi'$ is a Morse function,
then we can write explicit asymptotic expansions in
$\tau$ with coefficients from information about germs
of $\phi$ and $u$ at the critical points of $\phi'$.
Such asymptotic expansions are the basic in
determining
the asymptotic behavior of symbols under
product, adjoint, change of variables.
So they are the basic supporting
techniques in pseudo-differential calculus.
See \cite{Ho1}, chapter 7 for details.

\subsubsection{General phase function}

The phase function we are interested in has critical
points contained in
\begin{equation}
  \label{eq:cr:set}
  \{(x,\xi,g)\in T^*_1 M\x G: gx=x \text{ and } (dg)^*_x\xi=\xi\}.
\end{equation}
which at least contains the subset
\begin{equation}
  \label{eq:cr:sete}
  T^*_1 M\x \{e\}.
\end{equation}
In the above $T^*_1M$ means the quotient $(T^*M\mz)/\R_+$.
So the phase function is never a Morse function.

By recent results of \cite{Illman} and \cite{Kutz},
related to a generalization of Hilbert's seventh problem
to group action case, we may assume that the action
of $G$ is analytic on a compatible analytic atlas of $M$.
This atlas exists as shown in \cite{Illman},
and the uniqueness is proven in \cite{Kutz}.
The phase function is analytic under this atlas.

Qualitative results were introduced
for any analytic phase function,
using resolution of singularities
originated from Atiyah \cite{Atiyah:res}.
The asymptotic expansion formula has been
given by many people including Bernstein,
Malgrange, etc.. There is an excellent source for
of backgrounds, applications, explanations
of the result we use here.
It is in part II (page 169--268)
of volume II of the two volume monographs
by Arnold et al \cite{Arnold1},
\cite{Arnold2}. Although it is too
long a topic to elaborate on, all we need
is one technical theorem by Malgrange.

\subsubsection{Malgrange's theorem}

In \cite{Arnold2},  the following
theorem is derived and illustrated
with examples using Newton polyhedra
of the singularities.

\begin{theorem}\label{thm:osc}
Let $\phi$ be a real valued nonzero
analytic function on $\R^N$.
For $u\in C^\infty_c(\R^N)$ (real or complex
valued) a test function,
let $I(\tau)$ be the oscillatory integral
\begin{equation}
I_{\phi,u}(\tau)
=\int_{\R^N} e^{i\tau \phi(x)} u(x) dvol(x).
\end{equation}
Then for $\tau\to \infty$
\begin{equation}
I(\tau) = \sum\limits_{\alpha,p,q}
c_{\alpha,p,q}(u)\tau^{\alpha-p} (ln\,\tau)^q,
\end{equation}
where $\alpha\le 0$ runs through a finite set of
rational numbers, $p, q\in\N$ and $0\le q< n$.
Moreover $c_{\alpha,p,q}$ are all distributions
with support inside
\begin{equation}
S_\phi=\{x \in \R^N: d\phi(x)=0\},
\end{equation}
and with finite orders not exceeding $N$.
\end{theorem}

\subsection{Noncommutative residue}

\subsubsection{Definition}
In Connes-Moscovici local index formula, we need
to study the trace of $AP^{-z}$,  where $A\in\cPsDO_G$
and $P\in\cPsDO$ is a second-order, self-adjoint,
positive, elliptic pseudo-differential operator.
In fact $P$ is also $G$-invariant but we do not
use this assumption now. For convenience,
$P$ is assumed to be a second order, since
we will have any complex power of $P$ and
can easily change the scale of $z$.
We choose to work on an elliptic operator instead of
a transversally elliptic one 
because with all the other assumptions
we can add to $P$ a part $K$ which is transversally
smoothing and make $P+K$ is elliptic
(see section \ref{sec:red}). It will be shown that
such a $K$ has negligible contribution to the asymptotic
expansions we are interested in.

In this section we discuss the noncommutative residues
defined by such a $P$. In fact, to allow some room
for perturbation of $P$ relax the assumption on the
spectrum of $P$. In stead of asking 
\begin{equation}
  Spec(P)\subset (\epsilon,\infty) , 
\end{equation}
we allow an open neighborhood of it. This neighborhood
can be assumed to be a sector, as large as to miss
the negative half of the real axis. 

In summary we assume $P\in \cPsDO(M;E)$ is a second
order, elliptic, with 
\begin{equation}
  Spec(P)\subset \C\mz
  \sm\{\text{a conic neighborhood of }
  \R_-\} .
\end{equation}
Since $P$ is invertible and $P^{-1}$ is compact,
$Spec(P)$ is a discrete set.

For $\lambda$ on a contour along the
negative real axis $C$, the resolvant
$(P-\lambda)^{-1}$ exists and is bounded.
Let $\ln{z}$ denote the unique branch of 
the multiple valued logarithm function
defined on $\C\mz\sm\R_-$ such that the
imaginary part of it $arg(z)=0$ for
$z\in R_+$.
Then $\ln{(-\lambda)}$ is defined on
$\C\mz\backslash\R_+$ and with
imaginary part $arg(\ln(-\lambda))= 0$
on the negative half of the real axis.
Seeley (\cite{Seeley}, \cite{GS95}, \cite{Shubin})
showed that all the complex powers $P^z$ of $P$ are
pseudo-differential.

The function
\begin{equation}
  \label{eq:tr:conv}
  z\mapsto Trace(AP^{-z})
\end{equation}
is well defined and analytic on $Re(z)>\dim{M}$.

\begin{definition}
  For all $A\in\cPsDO^\infty(M;E)$, 
  assuming that $Trace(AP^{-z})$
    extends meromorphically the whole plane,
    we define
  \begin{equation}
    \label{eq:def:tau}
    \tau^P_k (A)= Res_{z=0} (z^k Trace(AP^{-z})).
  \end{equation}
\end{definition}

\subsubsection{Equivalent definitions}

\begin{lemma}\label{lem:itrans}
Suppose $N$ is a fixed positive integer,
$A\in\cPsDO^\infty(M)\rx G$,
$s_p$ are strictly decreasing sequences in $\R$.
Suppose for some fixed positive integer $K$
(typically $K$ is chosen to be $(\dim{M}+ord(A))/2$),
$A(P-\lambda)^{-K}$ is trace class.
Then the following three statements are equivalent:

(A) 
For $\lambda\to\infty$ on any ray of
$\C\backslash\R_+$,
\begin{equation}
\begin{split}
Tr(A(P-\lambda)^{-K})
\sim&
\sum_{l\ge 0}
\left(\sum_{q=0}^{N}
b_{p,q} (\ln{(-\lambda)})^q
\right) (-\lambda)^{-K+s_p}.
\end{split}
\end{equation}

(B) When $t\to 0+$
\begin{equation}
Tr(Ae^{-tP})
\sim
\sum_{p\ge 0}
\left(\sum_{q=0}^{N}
b'_{p,q} (\ln{t})^q\right) t^{-s_p}.
\end{equation}

(C) The zeta function $\Gamma(z)Tr(AP^{-z})$,
which converges for $Re(z)>K$,
extends to a meromorphic function on $\C$,
with the poles described as follows: 
up to a holomorphic function,
\begin{equation}
\Gamma(z)Tr(AP^{-z})
\sim
\sum_{p\ge 0}
\left(\sum_{q=0}^{N}
\frac{b''_{p,q}}{(z-s_p)^{q+1}}\right).
\end{equation}
Moreover, for any $p$, the coefficient subsets
$\{b_{p,q}: q=0, \ldots, N\}$,
$\{b'_{p,q}: q=0, \ldots, N\}$,
$\{b''_{p,q}: q=0, \ldots, N\}$ 
determine each other.
\end{lemma}


\begin{proof}
The proof is based on the following
straightforward integral transformation formulas 
(\cite{Gsch}):
\begin{equation}
\label{eq:mellin}
\begin{split}
P^{-s} =& \frac{1}{(s-1)\cdots(s-k)}
\frac{i}{2\pi}\int_{C_1} 
\lambda^{k-s}
\pd^k_\lambda(P-\lambda)^{-1}d\lambda\\
=&\frac{1}{\Gamma(s)}
\int_0^\infty t^{s-1} e^{-tP} dt,
\end{split}
\end{equation}
and
\begin{equation}
\label{eq:invmellin}
\begin{split}
e^{-tP} =& t^{-k}
\frac{i}{2\pi}\int_{C_1} 
e^{-t\lambda}
\pd^k_\lambda(P-\lambda)^{-1}d\lambda\\
=&\frac{1}{2\pi i}
\int_{Re(s)=c} t^{-s} \Gamma(s) P^{-s} ds,
\end{split}
\end{equation}
where $C_r, r\in (0,1]$ is the counterclockwise contour  
\begin{equation}
\begin{split}
C_r=&\{z\in\C: arg(z)=\pm(\pi/2-\delta), \abs{z}\ge r\}\\
&\cup\{z\in\C: arg(z)\in [-\pi/2+\delta,\pi/2-\delta],\abs{z}=r\}
=rC_1
\end{split}
\end{equation}
for some fixed small $\delta>0$.
Equations (\ref{eq:mellin}) (\ref{eq:mellin}) 
are the Mellin transform,
the inverse Mellin transform and
contour integrals that can be
converted to a slightly tilted version of
Mellin transform.
Including the appearance of higher logarithmic
terms, the argument has been given in \cite{Shubin},
see also \cite{Kassel:res}.
Here we only show the relation between the coefficients.
Note that the contribution from a different contour
only produce a smoothing operator.

To show (A) $\Rightarrow$ (B): 
we have
\begin{equation*}
A(P-\lambda)^{-s}
\in\cPsDO^{-K+ord(A)}(E,E)\rx G
\end{equation*}
is trace class by the choice of $K$.
By
\begin{equation}
\pd^k_\lambda(P-\lambda)^{-1}
=k!(P-\lambda)^{-k-1},
\end{equation}
we get
\begin{equation}
Ae^{-tP} =(N-1)! t^{-N+1}
\frac{i}{2\pi}\int_{C_1} 
e^{-t\lambda} A(P-\lambda)^{-N}d\lambda.
\end{equation}

Direct calculation shows
\begin{equation}
t^{-N+1}\frac{i}{2\pi}
\int_{C_1} e^{-t\lambda} (-\lambda)^{r_j-N}d\lambda
=\frac{t^{-r_j}}{\Gamma(-r_j)}
\end{equation}
and similarly
\begin{equation}
\begin{split}
t^{-N+1} &\frac{i}{2\pi}\int_{C_1}  e^{-t\lambda}
(-\lambda)^{s_p-N} (\ln{(-\lambda)})^p d\lambda\\
&= t^{-s_p} \frac{i}{2\pi}\int_{C_t} e^{-\lambda}
 (-\lambda)^{s_p} (\ln{(-\lambda)} -\ln{t})^q d\lambda\\
&\sim \sum_{i=0}^{q} (-1)^i \frac{q!}{i!(q-i)!} t^{-s_p} (\ln{t})^i
\frac{i}{2\pi}\int_{C_1} e^{-\lambda}
 (-\lambda)^{s_p} (\ln{(-\lambda)})^{q-i} d\lambda\\
&\sim \sum_{i=0}^{q} (-1)^i \frac{q!}{i!(q-i)!} t^{-s_p} (\ln{t})^i
(\frac{1}{\Gamma})^{(q-i)}(-s_p),
\end{split}
\end{equation}
where $(\frac{1}{\Gamma})^{(q-i)}$ is the $(q-i)$-$th$
derivative of the entire function $1/\Gamma$.
By the first equality of (\ref{eq:invmellin}),
and the fact that in this integral transform 
asymptotic expansions in $-\lambda$ correspond to
asymptotic expansions in $t$.

To see the relation between coefficient coefficients,
we look at the contribution of a single term,
\begin{equation}
b''_{p,q}=\sum_{i=0}^{q}
(-1)^i
\frac{q!}{i!(q-i)!} t^{-s_p}
(N-1)!
b_{p,i}
(\frac{1}{\Gamma})^{(q-i)}(-s_p).
\end{equation}
The above gives $b''_{p,q}$ recursively.

The rest of details can be found in,
for example, \cite{Shubin} and \cite{Kassel:res}.
\end{proof}

\subsection{Asymptotic trace formula for weakly polyhomogeneous
  parametrized pseudo-differential operators}

\subsubsection{Weakly parametric pseudo-differential operators}

For a self-adjoint, elliptic pseudo-differential operator $P$
of positive order, the analysis
of the resolvent $(P-\lambda)^{-1}$ for $\lambda$
in a certain sector $\Gamma$ of $\C\backslash\R_+$
plays a very important role in the study of
noncommutative residues.
In this section we recall some basic definition
and results by Grubb and Seeley \cite{GS95}.

Let $\Gamma$ be a sector of $\C\backslash\R_+$, of the form
\begin{equation}
\{\lambda=re^{i\theta} | r>0, \theta\in I\subset (0,2\pi)\},
\end{equation}
where $I$ is a subinterval, let $\Gamma^o$ be its
image under conjugacy which has the same form.
We are mainly interested in the case
when $\Gamma$ is near the negative real axis
but some discussion apply to any sector on $\C\mz$.

For simplicity we start with symbols in $M=\R^n$
and with scalar values. What matters for us is
the behavior of $r=\abs\mu\to\infty$
for $\mu\in\Gamma^o$,  or equivalently,
the behavior of $z=1/\mu\to 0$
for $z\in\Gamma$, which we use more often. 

\begin{definition}
The {\it weakly parametric} symbol space
$S^{m,0}(\R^n\x\R^n, \Gamma)$
consists of functions $p(x,\xi,\mu)$ that are
holomorphic in $\mu=1/z\in \Gamma^o$ for
\begin{equation}
\abs{\xi,\mu}=(\abs{\xi}^2+\abs{\mu}^2)^{1/2}
\ge \epsilon
\end{equation}
for some $\epsilon$ and satisfy, for all $j\in\N$,
$1/z\in \Gamma$,
\begin{equation}\label{GS11}
\pd^j_z p(\cdot,\cdot,1/z) \in S^{m+j}(\R^n\x\R^n)
\end{equation}
locally uniformly in $\abs{z}\le1$.
Moreover, for any $d\in\C$, define
\begin{equation}
S^{m,d}=\mu^dS^{m,0}.
\end{equation}
That is, $p(\cdot,\cdot,1/z)\in S^{m,d}$
if $z^dp(\cdot,\cdot,1/z)$ satisfies
(\ref{GS11}).
\end{definition}

$S^{m,d}$ are \Fr spaces with semi-norms
implied in the definition.

The following properties are immediate:

1. As symbols constant on $\Gamma$, $S^m\subset S^{m,0}$.

2. If $m<0$, then $f(x,\xi)\in S^m$ implies
$p(x,\xi,\mu)=f(x,\xi/\mu)\in S^{0,0}$.

3. $S^{m,d}\subset S^{m',d'}$ for $m\le m'$, $d'-d\in\N.$

4. The following maps are continuous
\begin{equation}
\begin{split}
&\pd^\alpha_\xi:  S^{m,d}\to S^{m-\abs\alpha,d}\\
& \pd^\beta_x:  S^{m,d}\to S^{m,d}\\
& z^k:  S^{m,d}\to S^{m,d-k}\\
& \pd^j_z:  S^{m,0}\to S^{m,0}\\
& \pd_z:  S^{m,0}\to S^{m+1,0}+S^{m,d+1}\; d\not=0
\end{split}
\end{equation}

5. $S^{m,d} \cdot S^{m',d'} \subset S^{m+m', d+d'}.$

As usual, let
\begin{equation}
S^{\infty,d}=\bigcup\limits_{m\in\R}S^{m,d},
\;\;\;\;
S^{-\infty,d}=\bigcap\limits_{m\in\R}S^{m,d}.
\end{equation}

\begin{definition}
Let $p_j$, $j\in\N$ be a sequence of
symbols in $S^{m_j,d}$, where $m_j\downarrow 0$.
Then an asymptotic expansion 
in $S^{m_0,d}$ (or $S^{\infty,d}$)
for $p\in S^{m_0,d}$
\begin{equation}
p\sim \sum_{j\in\N} p_j 
\end{equation}
means that for any $N\in\N$, 
\begin{equation}
p-\sum_{j=0}^N p_j\in S^{m_{N+1},d}.
\end{equation}
\end{definition}


The following is a generalization of
{\it classical polyhomogeneous} symbols,
which has asymptotic expansion of
homogeneous in $\xi$ of integer-degree
(see, for example, \cite{Shubin}).
Pseudo-differential operators with classical 
polyhomogeneous symbols are called
classical pseudo-differential operators.

\begin{definition}
$p\in S^{\infty,d}(\R^n\x\R^n,\Gamma)$ is called
{\it weakly polyhomogeneous} if there exists a sequence
of symbols $p_j\in S^{m_j-d,d}$, homogeneous in
$(\xi,\mu)$ for $\abs{\xi}\ge 1$ of degree
$m_j \downarrow -\infty$, such that
$p\sim \sum p_j\in S^{\infty,d}$.
\end{definition}
  
\begin{theorem}\label{gs112}
For $p\in S^{m,d}(\R^n\x\R^n,\Gamma)$,
the limits
\begin{equation}
p_{(d,k)}(x,\xi)
=\lim_{z\to 0}
\pd^k_z(z^d p(x,\xi,1/z)) \in S^{m-k}(\R^n\x\R^n)
\end{equation}
exists and for any $N\in\N$,
\begin{equation}
p(x,\xi,\mu) \sim \sum_{k=0}^N \mu^{d-k} p_{(d,k)}
\in S^{m+N+1,d-N-1}(\R^n\x\R^n,\Gamma).
\end{equation}
\end{theorem}



\begin{proposition}\label{GS:smooth}
For a parameterized pseudo-differential operator
$P(\lambda)$ on $\R^n$ with symbol
$p(x,\xi,\lambda)\in S^{-\infty,d}$, its kernel
has an expansion
\begin{equation}
K(x,y,\lambda)\sim
\sum_{k\in\N}
\lambda^{d-k}K_k(x,y),
\end{equation}
with $K_k$ smooth on $\R^n\x\R^n$, and
\begin{equation}
K-\sum_{k=0}^N K_k(x,y)\mu^{d-k}
\in C^\infty(\R^{2n}\x\Gamma),
\end{equation}
holomorphic in $\mu\in \Gamma^o$ for $\abs\mu\ge 1$.
\end{proposition}

In particular for smooth functions $\phi$ and $\psi$
with disjoint support, any $P(\lambda)$ with
symbol in $S^{m,d}$, $\phi P(\lambda) \psi$
is such an example.



\begin{definition}
A parameterized pseudo-differential operator $P(\mu)$
on $M$, $\mu\in\Gamma^o$ is in $\cPsDO^{m,d}(M)\x \Gamma$
if for any $\phi,\psi\in C^\infty_c(M)$, with support
in a common coordinate system $U$, 
$\phi P(\mu) \psi$ has a symbol in
$S^{m,d}(U\x\R^n\x\Gamma)$.
\end{definition}

So the noncommutative residues, are basically
the coefficients of the pure logarithmic grow
and its powers in the asymptotic expansion
of the resolvent.

\begin{definition} 
$\cPsDO^{-\infty,d}(M)\x\Gamma$,
consists of operators parameterized
in $\Gamma$ with kernel
$K(x,y,\mu)$, smooth in $(x,y)$, holomorphic
in $\mu$ for $\mu\in\Gamma^o$, $\abs\mu\ge 1$,
and has expansion
\begin{equation}
K(x,y,\mu)\sim \sum K_k(x,y,\mu)\mu^{d-k}
\end{equation}
with
\begin{equation}
\pd^\alpha_x\pd^\beta_\xi
(K-\sum_{k=0}^J K_k\mu^{d-k})
=O(\mu^{d-J-1})
\end{equation}
for $\abs\mu \to\infty$ locally uniformly. 
\end{definition}

All these concepts are of local natural, 
so the above symbolic calculus extends to
to the vector bundle version.
By a partition of unity subordinate to
a trivialization of the vector bundles,
we define
$S^{m,d}(M\x\Gamma; E, F)$
and $\cPsDO^{m,d}(M\x \Gamma; E,F)$.
in the usual way (see \cite{GS95} for details).

\begin{theorem}\label{GSShubin}
Let  $p(x,\xi)\in S^r(T^*M; E)$
be a weakly homogeneous symbol
of  positive order $m\in\N_+$ 
(i.e., $m$-homogeneous for $\abs\xi\ge 1$),
and assume $(p(x,\xi)+\mu^m)$ invertible
on a closed sector $\Gamma$,
then it is weakly polyhomogeneous on $\mu$,
with weakly parameterized symbols
\begin{equation}
(p(x,\xi)+\mu^m)^{-1}\in
S^{-m,0}(\Gamma)\cap S^{0,-m}(\Gamma).
\end{equation}
\end{theorem}

This theorem shows that Grubb and Seeley
generalized Seeley's (\cite{Seeley})
and Shubin's (\cite{Shubin}) methods
in spectral analysis for differential operators.
For example when a scalar symbol $p$
is strictly positive, there is a sector
$\Gamma$ near negative real axis such that
$\mu=(-\lambda)^{1/m}$ has a unique
root for $\lambda\in\Gamma$,
satisfying the invertibility condition.

\subsubsection{Asymptotic expansion of local trace formula}
For a family of pseudo-differential operators with
parameters $(g,\mu)$, uniformly in $G$,  and 
and weakly parametric in $\mu$
in a sector $\Gamma$, we find the
asymptotics the kernel and
trace of an operator in
$\cPsDO^{-n-1,d}(E,E; \Gamma)\rx G$.

In particular we are interested in traces of
operators of the form
\begin{equation}
  P (D-\lambda)^{-k}
\end{equation}
where $D$ is positive, positively ordered,
polyhomogeneous, elliptic pseudo-differential
operator and for what we concern $\lambda$ is
on the negative half of the real axis.

As explained in the introductory part of
section \ref{sec:32}, we may work on a
relatively compact neighborhood of $\R^n$.
By proposition \ref{GS:smooth}, the contribution
form a smoothing operator is already known,
we only consider operators modulo a smoothing
operator.

Also the scalar case is sufficient, in the
matrix case we take $tr_x$ and check the
invariance under transition functions. 

In $\R^n$, a pseudo-differential
operator can be written as the
Fourier integral form
\begin{equation}
Op(p)u(x)=\int\int
e^{i(x-y)\cdot\xi}
p(y,\xi) u(y) dyd\xi,
\end{equation}
and the distributional kernel is
the integral on $\R^n(\xi)$, which has possible
singularities only on the diagonal.

\begin{proposition}\label{key1}
Let $d\ge 0$, $U\subset\R^n$ be a relatively compact open chart
invariant under $G$ action,
\begin{equation}\begin{split}
p(x,\xi,g,\mu) &\sim \\
&\sum_{j=0}^\infty
p_j(x,\xi,g,\mu) \in S^{-n+d,d}(U\x\R^n,\Gamma,G)
\end{split}
\end{equation}
be a weakly polyhomogeneous in $(\xi,\mu)$, uniformly
in $G$, with $p_j$ homogeneous with degrees
$m_j\downarrow-\infty$.
Then $\int_G Op(p)\rho(g)dg$
has a continuous kernel
$K_p(x,y,\mu)$
with a diagonal expansion
\begin{equation}\label{eq:331}
K_p(x,x,\mu)
\sim
\sum_{j=0}^\infty c_j(x)\mu^{m_j+n}
+\sum_{\alpha,p,q} c'_{\alpha,p,q}(p)(x)
\mu^{d+\alpha-p}(\ln{\mu})^q
\end{equation}
for $\abs\mu\to\infty$, locally uniformly in $\Gamma$,
and for $\alpha, p,q$ described in Theorem \ref{thm:osc}.
All terms are decided by the symbol expect
$c'_{\alpha,p,0}$. The contribution to each $c_j$,
$c'_{\alpha,p,q}$ are from finitely many $p_j$'s.
\end{proposition}

The following proof is a adapted from the proof
of theorem 2.1 in Grubb-Seeley \cite{GS95}.

\begin{proof}
Without loss of generality we may assume $d=0$,
the general case follows by considering $\mu^d p$.

Let $r_J$ be remainder terms:
\begin{equation}
r_J=\sum_{j=0}^J p_j.
\end{equation}

$\int_G Op(r_J)\rho(g)\;dg$ has kernel
\begin{equation}
K_{r_J}(x,x,\mu)=
\int_{\R^n}\int_G e^{i(x-gx)\cdot\xi}
r_J(x,\xi,\mu) dg\;d\xi
\end{equation}
so by theorem \ref{gs112}, for any $N\in \N_+$ 
\begin{equation}
r_J(x,\xi,\mu) + \sum_{k=0}^{N} s_k(x,\xi)\mu^{-k}
+O(\abs{\xi}^{m_j-N}\mu^{-N}),
\end{equation}
where $s_k\in S^{m_j+k}$. For any $N$, taking large
enough $J$ will yield
\begin{equation}
K_{r_J}(x,x,\mu)=\sum_{k=0}^N c_{N,J,k}(x) \mu^{-k}
+O(\mu^{-N})
\end{equation}
so this will contribute like smoothing operators.

Now we discuss the contribution 
of each homogeneous symbols $p_j\in S^{m_j,0}$
which is homogeneous in $(\xi,\mu)$ of
degree $m_j$. As usual,
we split the integral into three parts:
\begin{equation}\begin{split}
K_{p_j}(x,x,\mu)
&=\int_{\R^n}\int_G e^{i(x-gx)\cdot\xi}
p_j(x,\xi,\mu) dg\;d\xi\\
&=\int_{\abs{\xi}\ge \abs{\mu}}
\int_G e^{i(x-gx)\cdot\xi}p_j(x,\xi,\mu) dg\;d\xi\\
&+\int_{\abs{\xi}\le 1}
\int_G e^{i(x-gx)\cdot\xi}p_j(x,\xi,\mu) dg\;d\xi\\
&+\int_{1\le \abs{\xi}\le \abs{\mu}}
\int_G e^{i(x-gx)\cdot\xi}p_j(x,\xi,\mu) dg\;d\xi.
\end{split}
\end{equation}

For the second and third integral, we use
theorem \ref{gs112} again,
\begin{equation}
p_j(x,\xi,\mu)
=\sum_{k=0}^M \mu^{-k}q_k(x,\xi)
+R_M(x,\xi,\mu),
\end{equation}
where
\begin{equation}
q_k(x,\xi)=\frac{1}{k!}
\pd^k_zp(x,\xi,1/z) \in S^{m_j+k}
\end{equation}
and homogeneous in $\xi$ for $\abs{\xi}\ge 1$
of degree $m_j+k$,
and 
\begin{equation}
R_M=O(\abs{\xi}^{m_j+M}\mu^{-M}).
\end{equation}

The second integral only contribute to
non-positive powers of $\mu$.

For the first integral, since $\abs\mu\ge 1$
is the only interested situation, $p_j$ is
homogeneous, so 
\begin{equation}\begin{split}
&\int_{\abs{\xi}\ge \abs{\mu}}
\int_G e^{i(x-gx)\cdot\xi}p_j(x,\xi,\mu) dg\;d\xi\\
&= \mu^{m_j+n} 
\int_{\abs{\xi}\ge 1}
\int_G e^{i \abs{\mu}(x-gx)\cdot\xi}
(\abs\mu^{-n-m_j}p_j)(x,\xi,\mu/\abs\mu) dg\;d\xi\\
&= \mu^{m_j+n} 
\int_1^\infty \int_{\abs{\xi'}= 1}
\int_G e^{i t \abs{\mu}(x-gx)\cdot\xi'}
(\abs\mu^{-n-m_j}p_j)(x,t \xi',\mu/\abs\mu) dg\;d\xi'\;dt.
\end{split}
\end{equation}
by the oscillatory integral theorem \ref{thm:osc} 
and then integrate over $t$ which are termwise
convergent, this gives terms (\ref{eq:331})
directly in the estimate.

The third integral break down into integrals of
for those of $q_k$, and $R_M$ which are all
homogeneous,
\begin{equation}\begin{split}
&\mu^{-k}\int_{1\le\abs\xi\le\abs\mu}\int_G
e^{i(x-gx)\cdot\xi} q_k(x,\xi) dg\;d\xi\\
&=\mu^{-k}\int_{1}^{\abs\mu} \tau^{m_j+k+n-1}
\int_{S^{n-1}}\int_G
e^{i\tau (x-gx)\cdot\xi'} q_k(x,\xi') dg\;d\xi'\;d\tau
\end{split}
\end{equation}
using polar coordinates in $\xi$-plane,
we get another part that contribute to
these logarithm expansions.

For $R_M$, which can be extended to homogeneous
in $\xi$ for all $\xi\not=0$, with a difference
of the second integral type, 
choose large enough $M>-n-m_j$, so that $r_M$
gives terms just as $q_k$, except also 
for non-positive powers of $\mu$.

The conclusion follows from the combination
of three integrals.

\end{proof}

\begin{proposition}\label{prop:main1}
Let $A\in \cPsDO^k(E,E)\rx G$, $P$
be a second order weakly polyhomogeneous,
self-adjoint and positive elliptic
pseudo-differential operator on $E$,
and let be $\Gamma$ a sector near
negative real axis.
Then for $Re(z)>(k+n)/2$, $AP^{-s}$ is
trace class,
$Trace(AP^{-s})$ is analytic in $z$ on
the half plane $Re(z)>(k+n)/2$, extending to 
a meromorphic function on $\C$ and up to an entire
function in $z\in\C$,
\begin{equation}
\Gamma(z)Tr(AP^{-z})
\sim \sum_{j\ge 0} 
\frac{\tilde{c}_j}{z+{\frac{j-k-n}{2}}}
+\sum_{l\ge 0}
\left(\sum_{p=0}^{\dim{G}+n-1}
\frac{\tilde{c}'_{p,q}}{(z-\frac{k+n}{2}+\frac{l}{q})^{p+1}}
\right),
\end{equation}
where all the coefficients except those
contributed also by the first sum
are determined by the symbol.
\end{proposition}

\begin{proof}
By theorem \ref{GSShubin}, $A(P+\mu^2)$
is weakly polyhomogeneous, so we apply
proposition \ref{key1}, for $-\lambda=\mu^2$.
Apply the integral transformation in lemma
\ref{lem:itrans}.
If $A'$ is smoothing,
then by corollary \ref{GS:smooth}
and proposition \ref{prop:main1}
\begin{equation}
Tr(A'(P-\lambda)^{-N})
\sim \sum_{j\ge 0}
c_j (-\lambda)^{-N-j/2}
\end{equation}
which does not contribute to the
pure logarithm and their power terms.
That is the generalized noncommutative
residues are locally computable.

\end{proof}

\subsection{Transversal residue formula}
\label{sec:red}

\begin{definition}
  For a pair of conic neighborhood $V_1$, $V_2$ 
  of $T^*_G M$ with $\bar{V_1} \subset V_2$,
  let 
  \begin{equation}
    \label{eq:atadd}
    \cPsDO^m_+(E, V_1, V_2)
  \end{equation}
  be those $P\in \cPsDO^m(E)$ such that
  \begin{itemize}
  \item[1.] $P$ is positive;
  \item[2.] $Ess(P) \cap V_1 =\ep ;$
  \item[3.] $P$ is $G$ invariant;
  \item[4.] If $\xi\in T^*M_x\mz \backslash V_2$
    and $\abs{\xi}>1$, then the symbol
    $\sigma(P)(x, \xi)$
    of $P$ is positive definite.
  \end{itemize}
\end{definition}

\begin{lemma}
  If $V_1$ and $V_2$ are small enough, 
  $\cPsDO^m_+(E, V_1, V_2)$ is nonempty.
\end{lemma}

\begin{proof}
For any $Y\in\g$, let
\begin{equation}
  Y_E: L^2(E)\to L^2(E)
\end{equation}
defined by the infinitesimal action of $\rho(exp(tY))$ on $E$,
as $t\to 0$. For any $Y\in \g$,
$Y_E$ is  a first order differential operator on $E$.

For an orthonormal basis $Y_i$, $i=1, \ldots, \dim{\g}$,
of $\g$, let $W_A$ be the differential operator 
introduced by Atiyah \cite{Atiyah:TEO}:
\begin{equation}
W_A= 1 - \sum_{i=1}^{\dim{\g}} Y_{i,E}^2.
\end{equation}
So the principal symbol of $W_A$ is positive
definite on $T^*M\mz\backslash T^*_GM$.
Atiyah (\cite{Atiyah:TEO}) used $W_A$  
to prove the existence of the distributional index.
Furthermore we can modify $W_A$ so that it achieves
the same but also transversally smoothing.

Choose a pair of small conic neighborhood $V_1, V_2$ of $T^*_G M$.
For this $V$ let $\chi$ and $P_\chi$ be the function and
pseudo-differential introduced in example \ref{eg:101}.
We assume $P_\chi$ is $G$-invariant, by averaging if
necessary. $W'_A=P_\chi W_A (P_\chi)*$ has essential support
disjoint with $T^*_G M$, and when $V_1$ and $V_2$ are small
enough $Q+W'_A$ is positive and elliptic.
In particular $W'_A\in \sK_G$. 
\end{proof}

Now take a $G$-invariant $Q\in \cPsDO^2(M)$
and assume it is positively ordered, transversally elliptic
and positive. It follows that the principal
symbol of $Q$ is positive definite on $T^*_G$.

\begin{proposition}
  Let $Q$ be as above.
  Then for any $W\in \cPsDO^2_+(E,V_1,V_2)$,
  $Q+W$ is elliptic, in addition to
  the above conditions.
  Moreover for all
  $k=0,1, \ldots, \dim{M}+\dim{G}-1$
  the residues $\tau_k^{Q+W}$ are
  independent of the choice of $W\in \cPsDO^m_+(E,V_1,V_2)$.
\end{proposition}

\begin{proof}

  When we add $W_A$ to $Q$ we get an operator $Q+W_A$
  which is $G$ invariant and elliptic.
  For $K_1$ and $K_2 \in \cPsDO^m_+(E,V_1,V_2)$, we have
  \begin{equation}
    (Q+K_1-\lambda)^{-1} -(Q+K_2-\lambda)^{-1} 
    =
    (Q-\lambda)^{-1} (K_1-K_2)(Q-\lambda)^{-1} .
  \end{equation}

  But the symbol right hand side is in $S^{-\infty,-2}$,
  so it has no contribution to the residues.
\end{proof}

\begin{definition}
  For $A\in \cPsDO^\infty(E)\rx G$, $Q$ as
  above, we define
  \begin{equation}
    \tau^{Q}_k(A)=\tau^{Q+K}_k(A)
  \end{equation}
  for some $K\in \cPsDO_+(E,V_1,V_2)$ and
  sufficiently small conic neighborhoods
  $V_1$ and $V_2$.
\end{definition}

In fact, more careful examination of the
modified operators, we find that they are
close to an invariant operator on $T^*_G$
direction tensored with an smoothing
operator on the complimentary direction.
In fact, Br\"unning and Heintz \cite{BrHe}
used similar techniques
on such operators and showed the following
theorem. Let $m$ be the dimension of $M_0/G$.
Let $(r, V)$ be any irreducible
representation of $G$, and let $\chi_r$ be
the character of $r$:
\[
\chi_r(g) = Trace_V(r(g)).
\]

\begin{theorem}
  Let $P\in \cPsDO^{2k}(E)$ be a positive,
  positively ordered, $G$ invariant and
  transversally elliptic.
  Let $\nu(\lambda)$ be the
  multiplicity of $r$ in the eigenspace of $P$
  with eigenvalue $\lambda$:
  \begin{equation}
    N_r(t)= \sum_{\lambda<t} \nu(\lambda).
  \end{equation}
  Then $N_r(t)$ is finite and as $t\to\infty$
  \begin{equation}
    \label{eq:bh32}
    N_r(t)\sim \frac{t^{\frac{m}{2k}}}{m (2\pi)^m}
    \int_{M_0} \frac{1}{vol(G_x)}
    \int_{T^*_{G,1} M} Tr_{(E_x\otimes V^*)^{G_x}}
    (\sigma(\xi_1)^{\frac{-m}{2k}}
    \otimes id_{V^*}) d\xi_1 dx,
  \end{equation}
  where $G_x$ is the stabilizer at $x$ and
  $dx$ is the volume form. The error term
  is $O(t^{\frac{m-1}{2k}} \log{t})$.
\end{theorem}

Let $c$ be the coefficient of $t^{\frac{m}{2k}}$.

By an integral transform:
\begin{equation}
  \label{eq:transn}
  Trace(\rho(\pi^*(\chi_r)) P^s)
  =\int_0^\infty t^z dN_r(t)
\end{equation}
we get 
\begin{equation}
  Res_{s=-m} Trace(\rho(\pi^*(\chi_r)) P^s)=c
\end{equation}

\begin{theorem}
Under the above assumptions,
for any Dixmier trace $Tr_\omega$ we have
\begin{equation}
  \begin{split}
    &Tr_\omega(\rho(\pi^*(\chi_r)) Q^{-m})\\
    &=\frac{1}{m (2\pi)^m}
    \int_{M_0} \frac{1}{vol(G_x)}
    \int_{T^*_{G,1} M} Tr_{(E_x\otimes V^*)^{G_x}}
    (\sigma(\xi_1)^{\frac{-m}{2k}}
    \otimes id_{V^*}) d\xi_1 dx.
  \end{split}
\end{equation}
In particular, all $\tau^{Q}_k=0$ for $k\ge 1$.
\end{theorem}

\section{The index of a transversally elliptic operator}

Let $P\in\cPsDO^k(E,F)$,
$P: \Gamma(E)\to \Gamma(F)$,
be a transversally elliptic pseudo-differential
operator relative to $G$ action on $E$ and $F$. 
Let
\begin{equation}
\pi_P : L^2(E,M) \to ker(P)
\end{equation}
be the projection to the kernel
of $P$. In this section we
define the index as a distribution
on $G$.

\subsection{Definition of the index}
\begin{lemma}\label{lemma:proj}
For any $u\in L^2(E)$,
$WF(\pi_P(u))\subset char(P)$.
\end{lemma}
\begin{proof}
By the regularity theorem \ref{thm:regularity},
\begin{equation}
WF(\pi_P(u))\subset WF(P(\pi_P(u))\cup char(P).
\end{equation}
But $P(\pi_P(u))=(P\pi_P)(u)=0$ and so it has empty wave front
set.

\end{proof}

\begin{theorem}\label{thm:index}
Let $P$ be a transversally elliptic operator.
For any $\phi\in\A$, $\rho(\phi) \pi_P$
is a smoothing operator. In particular,
for any $f\in C_c^\infty(G)$, $\rho(f) \pi_P$
is trace class.
\end{theorem}

\begin{proof}
$\rho(f)=\rho(\pi_2^*(f))$
is a special case of $\rho(\phi)$.
By Lemma \ref{lemma:proj},
Lemma \ref{lemma:xwf},
Theorem \ref{thm:compo},for any $u\in L^2(E)$,
\begin{equation*}
WF(\rho(\phi)\pi_Pu)\subset WF'(\rho(\phi))\circ WF(\pi_Pu).
\end{equation*}
Since $P$ is transversally elliptic,
$(T^*_GM)\cap char(P)=\ep$, the above
composition is empty.
So $\rho(\phi)\pi_P$ is a smoothing operator,
in particular, trace class.
\end{proof}

Now we are ready to introduce the definition
of the index of a transversally elliptic operator.
The following is equivalent to Atiyah's original definition
\cite{Atiyah:TEO}. 
Atiyah pointed out in \cite{Atiyah:TEO} 
the wave front set approach by \Hor
(see, for example, \cite{NeZi})
will show that 
the index below also makes sense
as a distribution on the Lie group
even for non-compact Lie groups.

\begin{definition}\label{def:index}
The index of a transversally elliptic operator
$P$ is defined to be the distribution on $G$
such that for any $f\in C_c^{\infty}(G)$,
\begin{equation}
index^G(P)(f)
=Trace(\rho(f)\pi_P)-Trace(\rho(f)\pi_{P^*}).
\end{equation}
\end{definition}

We observe that it is possible to define the index without
the $G$-invariant condition ( see \cite{Kor:GTEO1} and
\cite{Kor:GTEO2}).
Since $P$ is $G$-invariant under our assumption,
the index is a central distribution on $G$.

Definition \ref{def:index}
naturally extends to a distribution on the
groupoid $\G$ induced by the action of $G$
on $M$.

\begin{definition}\label{def:lindex}
We define the local index density
of a transversally elliptic operator
$P$ as the distribution on $\G_1=M\x G$
such that for any $\phi\in \A=C^\infty_c(M\x G)$,
\begin{equation}
index^\G(P)(\phi)
=Trace(\rho(\phi)\pi_P)-Trace(\rho(\phi)\pi_{P^*}).
\end{equation}
\end{definition}

When $G$ is the trivial group, then $index^\G$
is a smooth function on $M$ -- the index density
as in the heat kernel proof of
the index theorem for classical operators --
and its integral on $M$ is
the index of the operator $P$.
Although the index is a topological invariant, 
the index density is not.

\subsection{$K$-homology of the algebra $\A$}

We recall some definitions. A $p$-summable
pre-Fredholm module over $\A$ is a pair
$(\sH, F)$ where

(1) $\sH=\sH_+\oplus \sH_-$ is a
$\Z_2$-graded Hilbert space with grading
$\epsilon=1_{\sH_+}\oplus (-1_{\sH_{-}})$,
which is a $\Z_2$-graded left $\A$-module,

(2) $F\in B(\sH)$, $F\epsilon=-\epsilon F$,
$\phi(F^2-1)$ is compact for any $\phi\in\A$,

(3) for any $\phi\in\A$, 
$[F, \phi]\in \sL^p(\sH)$, the $p$-Schatten ideal
of compact operators.

A $p$-summable pre-Fredholm $\A$-module is
called a {\it $p$-summable Fredholm $\A$-module}
if in addition $F^2=1$.

Let
\begin{equation}
P: \Gamma^\infty(M,E) \to \Gamma^\infty(M,F),
\end{equation}
be a $G$-invariant
transversally elliptic pseudo-differential operator
of order 0. $P$ has a transversal parametrix $Q$,
which can be assumed to be $G$-invariant as well
(by averaging on $G$).
Let $\sH=L^2(E)\oplus L^2(F)$
with grading $\epsilon=1_{L^2(E)}\oplus (-1_{L^2(F)})$,
it is a graded $\A$-module through the action
$\rho_E\oplus \rho_F$. Let
\begin{equation}
F=\twobytwo 0QP0
\end{equation}

\begin{theorem}
Let $P$ in $\cPsDO^0(M;E,F)$ be
$G$-invariant and transversally elliptic.
Then the pre-Fredholm module $(\sH, F)$ introduced as above is
$p$-summable for all $p>\dim{M}$.
\end{theorem}

\begin{proof}

Condition (1) is obvious. (2) follows
from pseudo-local property of pseudo-differential
operators and $G$-invariance of $P$ and $Q$.

For (3), we have 
\begin{equation}
[F, \phi]\in \cPsDO^{-1})G(E\oplus F, E\oplus F)
\end{equation}
and by the $*$-algebra properties of $\cPsDO^{-1}_G(E\oplus F, E\oplus F)$
\begin{equation}
\abs{[F, \phi]}^{n+1}\in \cPsDO^{-n-1}_G(E\oplus F, E\oplus F)
\subset \sL^{1}{\sH}.
\end{equation}

\end{proof}

We now recall the standard process
to transform a pre-Fredholm $\A$-module into a
Fredholm $\A$-module, preserving $p$-summability
(for details see \cite{Co:NDG}, Appendix II of part I).

Given a pre-Fredholm $\A$ module $(\sH, F)$,
let $\tilde{\sH}=\sH\hat{\otimes}\sC$ be the
graded tensor product of $\sH$ with a $1+1$
dimensional graded Hilbert space
$\sC=\sC_+\oplus \sC_-$, with $\sC_\pm=\C$.
Then $\tilde{\sH}_+=\sH_+\oplus \sH_-$
and $\tilde{\sH}_-=\sH_-\oplus \sH_+$.
The $\A$-module structure on $\tilde{\sH}$
is given by
\begin{equation}
\tilde{\rho}(\phi)(\xi\tilde\otimes \eta)
=(\rho(\phi)\xi)\hat\otimes (\twobytwo{1}000\eta),
\end{equation}
for any $\phi\in\A$, $\xi\in\sH$, $\eta\in\sC$
($\A$ acts only non-trivially on $\tilde\sH_+$ as in $\sH$).
Since $F$ is of odd order, $\epsilon F=-F\epsilon$,
$F$ is alway of the form
\begin{equation}
F=\twobytwo 0QP0,
\end{equation}
so we define
\begin{equation}
\tilde F=\twobytwo 0{\tilde Q}{\tilde P} 0
\end{equation}
where
\begin{equation}
\tilde P=\twobytwo P{1-PQ}{1-QP}{(QP-2)Q}
\;\;\;\;
\tilde Q=\twobytwo {(2-QP)Q}{1-QP}{1-PQ}{-P}.
\end{equation}

\begin{proposition}
For a pre-Fredholm $\A$-module $(\sH,F)$,
$(\tilde\sH,\tilde F)$ is a Fredholm $\A$-module,
and there exists a pre-Fredholm $\A$-module
$(\sH_0, 0)$ with zero $\A$-action such that

(1) $\tilde\sH=\sH\oplus \sH_0$,

(2) for any $\phi\in\A$, $\phi(\tilde F - F\oplus 0)$ is compact.

Moreover, if $(\sH,F)$
is $p$-summable, so is $(\tilde\sH,\tilde F)$.
\end{proposition}       

In \cite{Co:NDG} (Part I, section 1),
Connes showed that there is
a trace $\tau: \A\to \C$,
\begin{equation}
\tau(\phi)=\frac{1}{2}Trace(\epsilon F[F,\phi]) 
\end{equation}
which gives the index map $K_0(\A)\to \Z$,
\begin{equation}
index F^+_e=(\tau\otimes Trace)(e)
\end{equation}
for any projection $e$ ($e=e^*=e^2$)
in the finite matrix algebra $M_q(\tilde{\A})$
for arbitrary $q$.

\begin{definition} (Connes, \cite{Co:NDG})
The trace $\tau$ (denoted by
$char(\tilde H, \tilde F)$)
is called the Connes
character of the 1-summable Fredholm
$\A$-module $(\tilde \sH, \tilde F)$.
\end{definition}

Next we show that
the character $\tau$
of $(\tilde \sH, \tilde F)$
is just like the
local index density.
This fact allows us to get
transversal index formula
by computing the Connes character.
 
\begin{proposition}
\begin{equation}
index^G(P)=\pi_*char(\tilde H, \tilde F),
\end{equation}
where $(\pi)_*=(\pi_G)_*$ is the push-forward
by the projection $\pi=\pi_G: M\x G \to G$.
\end{proposition}

\begin{proof}
In fact, for a $G$-invariant
transversally elliptic pseudo-differential
operator $P$, let $G$ be the Green operator
on $L^2(F)$, i.e.,
\begin{equation}
GP=1-\pi_P,
\;\;\;
PG=1-\pi_{P^*}.
\end{equation}
So the local index density can be written as
\begin{equation}
index^G(f)=Trace(\rho(\pi^*f)(1-GP))-Trace(\rho(\pi^*f)(1-PG)).
\end{equation}
For a transversal parametrix $Q$ of $P$, let
\begin{equation}
\tau'(\rho(\pi^*f))=Trace((\rho(\pi^*f)(1-QP))-Trace(\rho(\pi^*f)(1-PQ))
\end{equation}
which is well defined by proposition \ref{prop:parametrix}.
First, we show that 
\begin{equation}
\pi_*\tau'=index^G(P).
\end{equation}
From
\begin{equation}
\pi_P=(1-QP)\pi_P,
\;\;\;
\pi_{P^*}=\pi_{P^*}(1-PQ),
\end{equation}
we have
\begin{equation}
\begin{split}
index^G(P)(f)
=&Trace(\rho(\pi^*f)\pi_P)-Trace(\rho(\pi^*f)\pi_{P^*})\\
=&Trace(\rho(\pi^*f)(1-QP)\pi_P)
-Trace(\rho(\pi^*f)\pi_{P^*}(1-PQ))\\
=&\tau'(\pi^*f)-Trace(\rho(\pi^*f)(1-QP)GP)\\
 &\phantom{\tau'(\pi^*f)}
 +Trace(\rho(\pi^*f)PG(1-PQ))\\
=&\tau'(\pi^*f).
\end{split}
\end{equation}
The last equality in the above equation
holds since the last two terms are
\begin{equation}
\pm Trace(\rho(\pi^*(f))G(P-PQP)).
\end{equation}

Repeating the above process we conclude that for any
integer $n$,
\begin{equation}\label{c1pqn}
index^G(P)(\pi^*f)
=Trace(\rho(\pi^*f)(1-QP)^n)-Trace(\rho(\pi^*f)(1-PQ)^n).
\end{equation}

Direct computation shows that:
\begin{equation}\label{c2pqn}
\begin{split}
\tau(\pi^*f)
&=\frac{1}{2}Trace(
\twobytwo {\rho_E(\pi^*f)\oplus 0_F}00{-\rho_F(\pi^*f)\oplus 0_E}\\
&- \twobytwo 0{\tilde{Q}}{\tilde{P}}0
\twobytwo {\rho_E(\pi^*f)\oplus 0_F}00{-\rho_F(\pi^*f)\oplus 0_E}
\twobytwo 0{\tilde{Q}}{\tilde{P}}0)\\
&=Trace(\rho_E(\pi^*f)(1-QP)^2)-Trace(\rho_F(\pi^*f)(1-PQ)^2)\\
&=\tau'(\pi^*f).
\end{split}
\end{equation}
The constant factor $1/2$ is not essential,
it depends on our particular choice of the way which
a pre-Fredholm is transformed into
a Fredholm module.

For $p$-summable ($p>1$) Fredholm module, using some $n>p$
in (\ref{c1pqn}) and (\ref{c2pqn}).
\end{proof}

\subsection{The spectral triple in transversally elliptic case}

In this section we construct a spectral triple associated to 
a transversally elliptic pseudo-differential
operator $P$. We discuss only the even case to simplify
our argument, the odd case is similar. 
First we may assume it is of order $1$
as a pseudo-differential operator.

because when necessary we may multiply
$D$ with an appropriate power of
$1+\Delta$, this operation is
an isomorphism of Sobolev
spaces, hence it does not alter the index
of $D$.

Now let $\sH=L^2E \oplus L^2F$,
$\epsilon=1_E \oplus (-1)_F$,
and $\A$ acts on $\sH$ by
$\rho_E\oplus \rho_F$ as before.
Let $D$ be an operator of odd grading,
\begin{equation}
D=\twobytwo{0}{D_-}{D_+}{0}.
\end{equation}
We are interested in the case when $D$ is symmetric,
which is equivalent to $D_-=(D_+)^*$.
We now show that if $D$ is symmetric then $D$
is essentially self-adjoint.

For a symmetric pseudo-differential operator
\begin{equation}
D: \Gamma^\infty_c(M, E\oplus F) \to \Gamma^\infty(M, E\oplus F),
\end{equation}
we have an extension
\begin{equation}
D': \Gamma^{-\infty}(M, E\oplus F) \to \Gamma^{-\infty}_c(M, E\oplus F).
\end{equation}

\begin{proposition}
If $D\in \cPsDO(E\oplus F)$
is transversally elliptic,
symmetric and has positive order,
then as an unbounded operator on $\sH$,
defined on the domain of smooth sections,
$D$ is essentially self-adjoint.
\end{proposition}
Kordyukov \cite{Kor:GTEO2}) proved a more general
statement.

\begin{proof}

To show $D$ is essentially self-adjoint
we need only to show that (see, for instance,
theorem 26.1 of \cite{Shubin})
\begin{equation}
  \label{eq:dom0}
Ker(D^*\pm iI) \subset Dom(\bar{D}).
\end{equation}

First we recall the proof (\cite{Shubin}) that
\begin{equation}
  \label{eq:dom1}
  Dom(D^*)=\{s\in\sH; D's\in\sH\}.
\end{equation}
We denote by $W$ the right hand side of (\ref{eq:dom1}).
For all $u\in \Gamma^\infty$, $s\in\sH$, by definition
\begin{equation}
  \pair{Du,s} =   \pair{u,D's}.
\end{equation}
This implies $W \subset Dom(D^*)$ and $D'|_W=D^*|_W$.
If $u\in Dom(D^*)$, then there is a $w\in\sH$ such
that for any $v\in \Gamma^\infty$,
\begin{equation}
  \pair{u, Dv}=
  \pair{D^*u,v} =
  \pair{w,v}.
\end{equation}
But this implies $D'u=w$ so $w\in\sH$. In other words,
$u\in W$.

Next we show that
\begin{equation}
S=\{s \in W \,| \forall\phi \in \A\;\rho (\phi)s\in \Gamma^{\infty}
\}
\end{equation}
is contained in $Dom(\bar{D})$.

Let $\{f_n: n\in\N\}$ be a sequence of bump
functions on $G$ converging to the delta
distribution $\delta$ at the identity on $G$.
For $s\in W$, we have $\rho(\pi^*(f_n))s\to s$ in $\sH$
and
\begin{equation}
  \rho(\pi^*(f_n))D's = D'\rho(\pi^*(f_n))s\to D's
\end{equation}
in any $\sH$.
Recall that (the graph of) $\bar{D}$ is
defined by the closure of its graph.
So we have $s\in Dom{\bar{D}}$.

Now it is clear to prove (\ref{eq:dom0}) we need only
to show $Ker(D^* \pm i)\subset S$. The proof
for the two cases are essentially the same.
If $(D^* + i)s=0$, then $D's\in\sH$ and $(D' + i)s=0$.
Since $D'$ is of positive order, $D'+i$ and $D'$
have the same principal symbol.
So there is an $G$-invariant $Q$
that is a transversal parametrix for $D+i$:
$\rho (\phi)[1-Q(D+ i)]$ is smoothing.
Thus
\begin{equation}
\rho (\phi)s= \rho (\phi)(1-Q(D+ i))s +0\in \Gamma^\infty,
\end{equation}
which says $s\in S$. 

\end{proof}

From now on we may assume $D$ is self-adjoint,
replacing $D$ by its closure when necessary.

\begin{lemma}
  \label{lem:Dinv3}
For any nonzero real number $\lambda$,
$\rho (\phi)(D-\lambda i)^{-1}$ and
$(D-\lambda i)^{-1}\rho (\phi)$
are compact operators for all $\phi\in \A$.
\end{lemma}

\begin{proof}
Let $Q_\lambda$ be a transversal parametrix for $D-\lambda i$.
As discussed above, it might not necessarily
have negative order, but it can be chosen so
by a cutoff on the symbol in a conic neighborhood
of $T^*_GM$. So
$1-Q_\lambda(D-\lambda i)=K$
and
$1-(D-\lambda i)Q_\lambda=K'$.
Apply the inverse to the right hand side of
the first parametrix formula,
we have
\begin{equation}
(D-\lambda i)^{-1} -Q_\lambda = K(D-\lambda i)^{-1}.
\end{equation}
$Q_\lambda$ is bounded since it is
a pseudo-differential operator of 
order zero. It suffices to show that
$\rho (\phi)K(D-\lambda i)^{-1}$ is smoothing.
$(D-\lambda i)^{-1}$ is a bounded operator
from $H^{s}$ to $H^{s+1}$ so its composition
with a smoothing operator is still smoothing.
\end{proof}

\begin{proposition}\label{pps}
$(\A,\sH,D)$ is a $\dim{M}^+$-summable spectral triple.
\end{proposition}
\begin{proof}
Since $D$ commutes with $\rho(g)$,
For any $\phi\in\A$, 
\begin{equation}
[D,\phi]=\int_G [D, \phi(x,g)]\rho(g)d\mu(g)
\end{equation}
so it is compact, as $[D, \phi(x,g)]$ is
a pseudo-differential operator of negative
order. Lemma \ref{lem:Dinv3}
shows $(\A,\sH,D)$ is a spectral triple.
For any $\phi\in\A$, and $\lambda\in\C\sm\R$,
\begin{equation}
\phi (\lambda-\abs{D})^{-1}
= \phi Q_\lambda + \phi K_\lambda (\lambda-\abs{D})^{-1}
\in \sL^{(dim(M), \infty)}
\end{equation}
since the transversal parametrix can be
chosen to be a pseudo-differential operator
of order $-1$, $\phi K_\lambda$ is trace class
and $(\lambda-\abs{D})^{-1}$ is bounded.
\end{proof}

\begin{lemma}
\label{k1}
When $\abs{D}$ has a scalar principal symbol,
the spectral triple $(\A,\sH, D)$ is regular.
\end{lemma}

\begin{proof}
We need to show
\begin{equation}
\A\cup[D,\A]\subset Dom^\infty(\delta),
\end{equation}
and in fact we will show that
\begin{equation}
\A\cup[D,\A]\subset\cPsDO^0(E,E)\rx G \subset Dom^\infty(\delta),
\end{equation}
and $\delta$ preserves $\cPsDO^0(E,E)\rx G$. 

Elements in $\A$ have scalar symbols, so
elements $[D,\A]\in \cPsDO^1(E,E)\rx G$ have
vanishing principal symbols, which implies
$[D,\A]\in \cPsDO^0(E,E)\rx G$.

In general, since $\abs{D}$ has scalar symbols
of degree 1, the commutator with any classical
pseudo-differential operator has order,
at most zero. And since $\abs{D}$ is $G$-invariant,
the composition with $\rho(g)$ has no effect.
\end{proof}

\subsection{Asymptotic spectral analysis for spectral triples}

With an extra parameter on
pseudo-differential operators,
that is, for families of
pseudo-differential operators,
the above properties still hold.

By lemma \ref{k1} we have a regular spectral
triple when $\abs{D}$ has scalar principal
symbols.
As we showed in section 5.1, we need
to study the poles of
\begin{equation}
\zeta_A(z)=Tr(A\abs{D}^{-2z})
\end{equation}
where $A\in\A_D\subset \cPsDO^0(E,E)\rx G$.
For the purpose of estimation,
we will relax the condition on $A$,
only assuming  $A\in\cPsDO^0(E,E)$.

\subsection{Connes-Chern character in the
  periodic cyclic cohomology of $\A$}

As we have seen in the previous discussion,
the computation of the index for transversally
elliptic pseudo-differential operators
amounts to the computation of the Connes
character of a finitely summable
Fredholm module.
The Connes character evolved into
its new version that takes values in cyclic
cohomology.

For a pre-$C^*$-algebra $A$,
and a $p$-summable Fredholm $A$-module
$(H,F)$, let
\begin{equation}
Tr'(T)=\frac{1}{2}Trace(\epsilon F[F,T])
\end{equation} 
as the Connes character ($\tau(\pi^*f)=Tr'(\pi^*f)$).
The Connes-Chern character
in the periodic cyclic cohomology
\begin{equation}ch^*(H,F)\in HP^*(A)\end{equation}
is defined to be
\begin{equation}
ch^*(H,F)(a^0, \ldots, a^n)
=(-1)^{n(n-1)/2}\Gamma(\frac{n}{2}+1)
Tr'\left(a^0 [F, a^1]\cdots[F,a^n]\right)
\end{equation}
for $n$ even and
\begin{equation}
ch^*(H,F)(a^0, \ldots, a^n)
=\sqrt{2i}(-1)^{n(n-1)/2}\Gamma(\frac{n}{2}+1)
Tr'\left(a^0 [F, a^1]\cdots[F,a^n]\right)
\end{equation}
for $n$ odd (Here $\Gamma$ is the
Gamma function).

As before, let $e\in M_q(\tilde{\A})$ be 
a projection (which is equivalent
to a finitely generated projective
module on $\A$ in the
$K$-theory $K_*(\A)$
for operator algebras
). There is a well understood
Chern character (see \cite{Co:book},\cite{Bla})
from $K_*(\A)$ to the periodic cyclic homology
of $\A$:
\begin{equation}
ch_*(e) \in HP_*(\A).
\end{equation}
The Connes character
can be viewed the dual of this
Chern character.
The Connes-Chern character gives the
index formula (for example, in the even case)
in the following fashion.
For an element $[e]\in K(A)$, then the
index of $F^+_e$ the twisted operator $F^+$
by the the projection $e$ is
\begin{equation}
Index F^+_e = \langle ch^*(H,F), ch_*(e)\rangle.
\end{equation}

In conclusion, to find the index of $P$
and its twisted versions, we may
compute the Connes-Chern character 
of a Fredholm module $F$ associated with it,
in the periodic cyclic cohomology.

\subsection{The Connes-Moscovici local index formula}
First, we briefly recall the definitions and results
from \cite{CM:local} (also see \cite{Co:book},\cite{GVF}).

Let $\A$ be a $*$-algebra, which
is a dense $*$-subalgebra of a pre-$C^*$
algebra $A$.

\begin{definition}
{\it An spectral triple}
is a triple $(\A, \sH, D)$ where

(1) $\sH=\sH_+\oplus\sH_-$
is a $\Z_2$-graded Hilbert space
and left $\A$-module, with
grading $\epsilon$;

(2) $D$ is an unbounded self-adjoint
operator on $\sH$ such that $D\epsilon=-\epsilon D$;

(3) for all $a\in\A$, $[D,a]\in B(\sH)$;

(4) for all $a\in \A$, $a(1+D^*D)^{-1}$ is compact.

{\it An odd spectral triple} over $A$
is similarly defined except
without grading and grading related conditions
in the above.
\end{definition}

Starting with a spectral triple, by the observation
in \cite{BaaJulg} by Baaj and Julg,
the following assignment
\begin{equation}
D \mapsto D(1+D^*D)^{-1/2} 
\end{equation}
determines a pre-Fredholm module.
In fact it is shown in \cite{BaaJulg}
that all $K$-homology classes
of a pre-$C^*$ algebra can be
obtained this way.  
We may switch to the computation
of the Connes-Chern character of
a spectral triple for the following
reasons:
the index of $D$ is preserved in the
Baaj-Julg assignment;
the Connes-Chern character is
an invariant of the $K$-homology
class; and when $D$ is
a pseudo-differential operator,
so is $D(1+D^*D)^{-1/2}$ whose symbol
can be computed in terms of
the symbol of $D$. 

Let $F$ and $\abs{D}$ be the elements
of the polar decomposition of $D$:
\begin{equation}
D=F\abs{D},
\end{equation}
where $F=sign{D}$ is unitary
and $\abs{D}=(D^2)^{1/2}$ is positive.
$(\sH,F)$ is a bounded
pre-Fredholm module, representing the
same class in $K$-homology
determined by the spectral triple.
For what we do, the invertibility of
$\abs{D}$ is not so important, a small
shift of its spectrum by a positive
number will not affect the asymptotic
behavior of the spectrum.
As we have showed, for our purpose,
$\abs{D}^{-1}$ only need to be well
defined up to a smoothing operator.
For example the Green operator works
fine:
\begin{equation}
\abs{D} G = 
G \abs{D} = 1 - \pi_{\ker{D}}.
\end{equation}

\begin{definition}
For some $p\ge 1$, 
a $p^+$-summable spectral triple is
a spectral triple $(\A,\sH,D)$
such that for any $\phi\in\A$
and $\lambda\in\C\sm\R$,
\begin{equation}
\phi(\lambda-\abs{D})^{-1}\in\sL^{p+}(\sH)=\sL^{(p,\infty)}(\sH),
\end{equation}
where for $p>1$, $\sL^{(p, \infty)}$ is the ideal
of $B(\sH)$ consisting of those
compact operators $T$ whose $n$-$th$
characteristic value
$\mu_n(\abs{T})=min\{\; \norm{\,\abs{T}\, |_{E^\perp}}; \dim{E}=n\}$
satisfies
\begin{equation}
\mu_n(\abs{T}) = O(n^{-1/p}).
\end{equation}
When $p=1$, this is a sufficient condition
for $T\in \sL^{(1, \infty)}$, but usually this is
satisfied.
\end{definition}

Let $\delta$ be the derivation
operator $ad(\abs{D})$:
\begin{equation}
\delta(T)=ad(\abs{D})(T)=[\abs{D}, T]
\end{equation}
defined on the bounded operators $B(\sH)$
and takes values as unbounded operators on $\sH$.
Let $Dom(\delta)\subset B(\sH)$ be the domain
of $\delta$; that is,
$A\in B(\sH)$ is in $Dom(\delta)$ if
and only if $[\abs{D},A]$ extends
to a bounded operator on $\sH$. And let
\begin{equation}
Dom^\infty(\delta)=\bigcap_{k\ge 1} Dom(\delta^k).
\end{equation}

\begin{definition}
A $p^+$-summable spectral triple is regular  if
\begin{equation}
\A\cup[D,\A]\subset Dom^\infty(\delta)
\end{equation}
\end{definition}

For a $n^+$-summable spectral triple $(A, H,D)$
which satisfies
\begin{equation}
\A\cup[D,\A]\subset Dom(\delta^2),
\end{equation}
Connes character formula says
the $n$-cocycle in Hochschild cohomology of $\A$ is,
for $a_i\in A$,
\begin{equation}\label{eq:char0}
\phi_\omega(a^0, \ldots, a^n)
=\lambda_n
Tr_\omega(\epsilon a^0[D,a^1]\cdots[D,a^n]\abs{D}^{-n}),
\end{equation}
where $\lambda_n$ is a universal constant,
and $Tr_\omega$ is the Dixmier trace (see
\cite{Co:book} IV.2. for details).
In the odd case the same is true with $\epsilon=id$.

Connes-Moscovici \cite{CM:local}
went further to solve the general
problem for local character formula
for spectral triples with
discrete dimension spectrums.
We now to introduce more notations.
Let $\A_D$ be a subspace of $B(\sH)$
generated by the following operations:
for any operators $A\in \A$,
\begin{equation*}
dA=[D,A],\;
\nabla(A)=[D^2,A],\;
A^{(k)}=\nabla^k(A)
\end{equation*}
are in $\A_D$, moreover
the operators
\begin{equation*}
P(a^0,a^1,\ldots,a^n)=
a^0 (da^1)^{(k_1)} \ldots (da^n)^{(k_n)}
\end{equation*}
are in $\A_D$,
where $a^0,\ldots, a^n \in \A$,
acting on $\sH$ by $\rho$.

To apply our result for the residues
the Connes-Moscovici theorems,
we first list the following criterions
in our specific spectral triple.

{\bf Criterions} We assume that
$(\A,\sH,D)$ is (1) $p^+$-summable for
some $p>1$, (2) regular,
(3) with discrete dimension spectrum,
(4) for $P\in\A_D$, the zeta
function
\begin{equation*}
\zeta_{P,D}(z)=Trace(P \abs{D}^{-2z})
\end{equation*}
is at least defined and analytic for
$Re(z)>k+p$. In this case 
\begin{equation}
\tau_q(P)=\tau^{\abs{D}}_q(P)=Res_{z=0} z^q \zeta_{P,D}.
\end{equation}
(5) Only finite many of $\tau_q$, $q=0, \ldots$
are nonzero in general.

We have seeing that our spectral $(\A,\sH,D)$ 
for the transversally elliptic case is regular,
$p^+$-summable, with discrete dimension
spectrum, and the poles of
$\zeta_{P,D}$ have multiplicities
not exceeding a  fixed number. Therefore we have
the following.

\begin{theorem} (Connes-Moscovici)
a) The following formula defines an
even cocycle in $(b,B)$ bicomplex of $\A$:
\begin{equation}
\begin{split}
\phi_0(a_0)&=\tau_{-1}(\gamma a_0)\\
\phi_{2m}(a_0,\ldots,a_{2m})
&=\sum_{k\in\N^{2m},q\ge 0}
c_{2m,k,q}\\
&\tau_q(\gamma a_0(da_1)^{(k_1)}\ldots(da_{2m})^{k_{2m}}
\abs{D}^{-2\abs{k}-2m})
\end{split}
\end{equation}
for $m>0$, where $c_{2m,k,q}$ are universal
constants given by
\begin{equation}\label{eq:evenconst}
c_{2m,k,q}=
\frac{(-1)^{\abs{k}}}{k!\tilde{k}!}\sigma_q(\abs{k}+m)
\end{equation}
where $k!=k_1!\ldots k_{2m}!$,
$\tilde{k}!=(k_1+1)(k_1+k_2+2)\ldots (k_1+\ldots+k_{2m}+2m)$,
and $\sigma_q(N)$ is the $q$-$th$ elementary polynomial
of the set $\{1,2,\ldots,N-1\}$.

(b) The cohomology class of the cocycle
$(\phi_{2m})_{m\ge 0}$ in $HC^{ev}(\A)$
coincides with the Connes-Chern character
$ch_*(\A,\sH,D)$.
\end{theorem}

In the odd case, the Connes-Chern
character is computed similarly. Suppose
we have an odd spectral triple still called
$(\A,\sH,D)$.

\begin{theorem} (Connes-Moscovici)
a) The following formula defines an
odd cocycle in $(b,B)$ bicomplex of $\A$:
\begin{equation}
\begin{split}
\phi_{2m+1}(a_0,\ldots,a_{2m+1})
&=\sqrt{2i}\sum_{k\in\N^{2m+1},q\ge 0}
c_{2m+1,k,q}\\
&\tau_q(a_0(da_1)^{(k_1)}\ldots(da_{2m+1})^{k_{2m+1}}
\abs{D}^{-2\abs{k}-2m-1})
\end{split}
\end{equation}
where $c_{2m+1,k,q}$ are universal
constants given by
\begin{equation}
c_{2m+1,k,q}=
\frac{(-1)^{\abs{k}}}{k!\tilde{k}!q!}\Gamma^{(q)}(\abs{k}+m+\frac12)
\end{equation}
where $\Gamma^{(q)}$ is the $q$-$th$ derivative of
the Gamma function.

(b) The cohomology class of the cocycle
$(\phi_{2m+1})_{m\ge 0}$ in $HC^{od}(\A)$
coincides with the Connes-Chern character
$ch_*(\A,\sH,D)$.
\end{theorem}

\begin{remark}
As a sample, where $G$ is the trivial group.
This is the commutative case since
$\A=C^\infty(M)$. $D$ is
an elliptic pseudo-differential operator
with scalar principal symbol
on a hermitian bundle $E$
over a smooth manifold $M$.
There exists such geometric operator
$D$, since we can choose a space-time
oriented pseudo-Riemannian manifold
with a compact Lie group acting
by isometries,  assume it is spin, then
the natural Dirac operator is obviously
not elliptic. But if there is no
isotropic (null) directions along the
orbit, then the Dirac operator is
transversally elliptic relative to
the group action (see \cite{Baum}).
The spectral triple is regular,
$\dim{M}^+$-summable, and with
discrete dimension spectrum.
Since the zeta functions involved
have at most simple poles, only $\tau_0$
is possibly nonzero. $\tau_0$ is up to
constant factors the Wodzicki
residue, which is computable
in terms of symbols of the operators
involved.
\end{remark}

\bibliographystyle{plain}

\begin{thebibliography}{10}

\bibitem{Arnold1}
V.~I. Arnol{\cprime}d, S.~M. Guse\u{\i}~n Zade, and A.~N. Varchenko.
\newblock {\em Singularities of differentiable maps. {V}ol. {I}}, volume~82 of
  {\em Monographs in Mathematics}.
\newblock Birkh\"auser Boston Inc., Boston, MA, 1985.
\newblock The classification of critical points, caustics and wave fronts,
  Translated from the Russian by Ian Porteous and Mark Reynolds.

\bibitem{Arnold2}
V.~I. Arnol{\cprime}d, S.~M. Guse\u{\i}~n Zade, and A.~N. Varchenko.
\newblock {\em Singularities of differentiable maps. {V}ol. {II}}, volume~83 of
  {\em Monographs in Mathematics}.
\newblock Birkh\"auser Boston Inc., Boston, MA, 1988.
\newblock Monodromy and asymptotics of integrals, Translated from the Russian
  by Hugh Porteous, Translation revised by the authors and James Montaldi.

\bibitem{Atiyah:res}
M.~F. Atiyah.
\newblock Resolution of singularities and division of distributions.
\newblock {\em Comm. Pure Appl. Math.}, 23:145--150, 1970.

\bibitem{Atiyah:TEO}
M.~F. Atiyah.
\newblock {\em Elliptic operators and compact groups}.
\newblock Springer-Verlag, Berlin, 1974.
\newblock Lecture Notes in Mathematics, Vol. 401.

\bibitem{AB:Lef1}
M.~F. Atiyah and R.~Bott.
\newblock A {L}efschetz fixed point formula for elliptic complexes. {I}.
\newblock {\em Ann. of Math. (2)}, 86:374--407, 1967.

\bibitem{AS1}
M.~F. Atiyah and I.~M. Singer.
\newblock The index of elliptic operators. {I}.
\newblock {\em Ann. of Math. (2)}, 87:484--530, 1968.

\bibitem{BaaJulg}
Saad Baaj and Pierre Julg.
\newblock Th\'eorie bivariante de {K}asparov et op\'erateurs non born\'es dans
  les ${C}\sp{\ast} $-modules hilbertiens.
\newblock {\em C. R. Acad. Sci. Paris S\'er. I Math.}, 296(21):875--878, 1983.

\bibitem{Baum}
Helga Baum.
\newblock The index of the pseudo-{R}iemannian {D}irac operator as a
  transversally elliptic operator.
\newblock {\em Ann. Global Anal. Geom.}, 1(2):11--20, 1983.

\bibitem{BV}
Nicole Berline and Mich{\`e}le Vergne.
\newblock The equivariant {C}hern character and index of ${G}$-invariant
  operators. {L}ectures at {C}{I}{M}{E}, {V}enise 1992.
\newblock In {\em $D$-modules, representation theory, and quantum groups
  (Venice, 1992)}, Lecture Notes in Math., 1565, pages 157--200. Springer,
  Berlin, 1993.

\bibitem{Bla}
Bruce Blackadar.
\newblock {\em ${K}$-theory for operator algebras}.
\newblock Cambridge University Press, Cambridge, second edition, 1998.

\bibitem{BrHe}
Jochen Br{\"u}ning and Ernst Heintze.
\newblock Representations of compact {L}ie groups and elliptic operators.
\newblock {\em Invent. Math.}, 50(2):169--203, 1978/79.

\bibitem{CM:local}
A.~Connes and H.~Moscovici.
\newblock The local index formula in noncommutative geometry.
\newblock {\em Geom. Funct. Anal.}, 5(2):174--243, 1995.

\bibitem{Co:NDG}
Alain Connes.
\newblock Noncommutative differential geometry.
\newblock {\em Inst. Hautes \'Etudes Sci. Publ. Math.}, 62:257--360, 1985.

\bibitem{Co:book}
Alain Connes.
\newblock {\em Noncommutative geometry}.
\newblock Academic Press Inc., San Diego, CA, 1994.

\bibitem{Cordes}
H.~O. Cordes.
\newblock {\em The technique of pseudodifferential operators}, volume 202 of
  {\em London Mathematical Society Lecture Note Series}.
\newblock Cambridge University Press, Cambridge, 1995.

\bibitem{Duistermaat:FIO}
J.~J. Duistermaat.
\newblock {\em Fourier integral operators}.
\newblock Birkh\"auser Boston Inc., Boston, MA, 1996.

\bibitem{GGK}
Israel Gohberg, Seymour Goldberg, and Nahum Krupnik.
\newblock {\em Traces and determinants of linear operators}, volume 116 of {\em
  Operator Theory: Advances and Applications}.
\newblock Birkh\"auser Verlag, Basel, 2000.

\bibitem{GVF}
Jos{\'e}~M. Gracia-Bond{\'\i}a, Joseph~C. V{\'a}rilly, and H{\'e}ctor Figueroa.
\newblock {\em Elements of noncommutative geometry}.
\newblock Birkh\"auser Boston Inc., Boston, MA, 2001.

\bibitem{Gsch}
Gerd Grubb and Elmar Schrobe.
\newblock Trace expansions and the noncommutative residue for manifolds with
  boundary.
\newblock {\em J. reine angew. Math.}, 536:167--207, 2001.

\bibitem{GS95}
Gerd Grubb and Robert~T. Seeley.
\newblock Weekly parametric pseudodifferential operators and
  {A}tiyah-{P}atodi-{S}inger boundary problems.
\newblock {\em Inventiones Mathematicae}, 121:481--529, 1995.

\bibitem{Ho1}
Lars H{\"o}rmander.
\newblock {\em The analysis of linear partial differential operators. {I}},
  volume 256 of {\em Grundlehren der Mathematischen Wissenschaften [Fundamental
  Principles of Mathematical Sciences]}.
\newblock Springer-Verlag, Berlin, 1983.
\newblock Distribution theory and Fourier analysis.

\bibitem{Ho3}
Lars H{\"o}rmander.
\newblock {\em The analysis of linear partial differential operators. {III}},
  volume 274 of {\em Grundlehren der Mathematischen Wissenschaften [Fundamental
  Principles of Mathematical Sciences]}.
\newblock Springer-Verlag, Berlin, 1994.
\newblock Pseudo-differential operators.

\bibitem{Illman}
S{\"o}ren Illman.
\newblock Every proper smooth action of a {L}ie group is equivalent to a real
  analytic action: a contribution to {H}ilbert's fifth problem.
\newblock In {\em Prospects in topology (Princeton, NJ, 1994)}, pages 189--220.
  Princeton Univ. Press, Princeton, NJ, 1995.

\bibitem{Julg}
Pierre Julg.
\newblock Induction holomorphe pour le produit crois\'e d'une ${C}\sp{\ast}
  $-alg\`ebre par un groupe de {L}ie compact.
\newblock {\em C. R. Acad. Sci. Paris S\'er. I Math.}, 294(5):193--196, 1982.

\bibitem{Kassel:res}
Christian Kassel.
\newblock Le r\'esidu non commutatif (d'apr\`es {M}.\ {W}odzicki).
\newblock {\em Ast\'erisque}, (177-178):Exp.\ No.\ 708, 199--229, 1989.
\newblock S\'eminaire Bourbaki, Vol.\ 1988/89.

\bibitem{Kor:GTEO1}
Yuri{\u\i}~A. Kordyukov.
\newblock Transversally elliptic operators on ${G}$-manifolds of bounded
  geometry.
\newblock {\em Russian J. Math. Phys.}, 2(2):175--198, 1994.

\bibitem{Kor:GTEO2}
Yuri{\u\i}~A. Kordyukov.
\newblock Transversally elliptic operators on a ${G}$-manifold of bounded
  geometry.
\newblock {\em Russian J. Math. Phys.}, 3(1):41--64, 1995.

\bibitem{Kutz}
Frank Kutzschebauch.
\newblock On the uniqueness of the analyticity of a proper ${G}$-action.
\newblock {\em Manuscripta Math.}, 90(1):17--22, 1996.

\bibitem{MRW:groupoidcstar}
Paul~S. Muhly, Jean~N. Renault, and Dana~P. Williams.
\newblock Equivalence and isomorphism for groupoid ${C}\sp \ast$-algebras.
\newblock {\em J. Operator Theory}, 17(1):3--22, 1987.

\bibitem{NeZi}
A.~Nestke and F.~Zickermann.
\newblock The index of transversally elliptic complexes.
\newblock In {\em Proceedings of the 13th winter school on abstract analysis
  (Srn\'\i, 1985)}, volume~9, pages 165--175 (1986), 1985.

\bibitem{BEP}
Bent~E. Petersen.
\newblock {\em Introduction to the {F}ourier transform and pseudodifferential
  operators}.
\newblock Pitman (Advanced Publishing Program), Boston, MA, 1983.

\bibitem{Schweitzer}
Larry~B. Schweitzer.
\newblock Dense {$m$}-convex {F}r\'echet subalgebras of operator algebra
  crossed products by {L}ie groups.
\newblock {\em Internat. J. Math.}, 4(4):601--673, 1993.

\bibitem{Seeley}
R.~T. Seeley.
\newblock Complex powers of an elliptic operator.
\newblock In {\em Singular Integrals (Proc. Sympos. Pure Math., Chicago, Ill.,
  1966)}, pages 288--307. Amer. Math. Soc., Providence, R.I., 1967.

\bibitem{Shubin}
M.~A. Shubin.
\newblock {\em Pseudodifferential Operators and Spectral Theory}.
\newblock Springer-Verlag, Berlin, 1980.

\bibitem{Smagin-Shubin}
S.~A. Smagin and M.~A. Shubin.
\newblock The zeta function of a transversally elliptic operator.
\newblock {\em Uspekhi Mat. Nauk}, 39(2(236)):187--188, 1984.

\bibitem{Stinespring}
W.~Forrest Stinespring.
\newblock A sufficient condition for an integral operator to have a trace.
\newblock {\em J. Reine Angew. Math.}, 200:200--207, 1958.

\bibitem{MTaylor}
Michael~E. Taylor.
\newblock {\em Pseudodifferential operators}, volume~34 of {\em Princeton
  Mathematical Series}.
\newblock Princeton University Press, Princeton, N.J., 1981.

\end{thebibliography}

\def\cprime{$'$} \def\cprime{$'$}

\noindent Address:\\
100 St. George St.\\
Department of Mathematics\\
University of Toronto\\
Toronto, Ontario Canada M5S 3G3

\begin{verbatim}
Email: xdhu@math.toronto.edu
\end{verbatim}

\end{document}